\documentclass[a4paper, reqno, 14pt]{amsart}

\usepackage[usenames,dvipsnames]{color}
\usepackage{amsthm,amsfonts,amssymb,amsmath,amsxtra}
\usepackage[all]{xy}
\SelectTips{cm}{}
\usepackage{xr-hyper}
\usepackage[colorlinks=
   citecolor=Black,
   linkcolor=Red,
   urlcolor=Blue]{hyperref}
\usepackage{verbatim}

\usepackage[margin=1.25in]{geometry}
\usepackage{mathrsfs}

\usepackage{diagbox}

\RequirePackage{xspace}
\RequirePackage{etoolbox}
\RequirePackage{varwidth}
\RequirePackage{enumitem}
\RequirePackage{tensor}
\RequirePackage{mathtools}
\RequirePackage{longtable}
\RequirePackage{multirow}

\usepackage{tikz}
\usepackage{tikz-cd}
\usepackage{csquotes}

\usepackage[url=false,doi=false,isbn=false,sorting=nyt,giveninits=true,maxbibnames=4]{biblatex}
\bibliography{references.bib}

\def\ge{\geqslant}
\def\le{\leqslant}
\def\a{\alpha}

\def\G{\Gamma}
\def\d{\delta}

\def\o{\omega}

\def\s{\sigma}

\def\k{\kappa}
\def\l{\lambda}

\def\i{^{-1}}

\def\<{\langle}
\def\>{\rangle}

\newcommand{\bJ}{\mathbf J}
\newcommand{\de}{{\mathrm{def}}}

\newcommand{\bG}{\mathbf G}

\newcommand{\BF}{\ensuremath{\mathbb {F}}\xspace}
\newcommand{{\BG}}{\ensuremath{\mathbb {G}}\xspace}

\newcommand{\BJ}{\ensuremath{\mathbb {J}}\xspace}
\newcommand{{\BK}}{\ensuremath{\mathbb {K}}\xspace}

\newcommand{\BS}{\ensuremath{\mathbb {S}}\xspace}

\newcommand{\CO}{\ensuremath{\mathcal {O}}\xspace}

\newcommand{\Ad}{{\mathrm{Ad}}}

\DeclareMathOperator{\Aut}{Aut}

\DeclareMathOperator{\diag}{diag}

\DeclareMathOperator{\Gal}{Gal}

\DeclareMathOperator{\Res}{Res}

\newcommand{\SL}{{\mathrm{SL}}}

\DeclareMathOperator{\LP}{LP}

\def\tW{\tilde W}

\def\doubleparenthesis#1{{(\!({#1})\!)}}
\def\doublebracket#1{{[\![{#1}]\!]}}

\newtheorem{theorem}{Theorem}

\theoremstyle{definition}

\newtheorem{example}[theorem]{Example}
\newtheorem*{example*}{Example}

\newtheorem*{function*}{Function}

\numberwithin{equation}{section}
\numberwithin{theorem}{section}

\setitemize[0]{leftmargin=*,itemsep=\the\smallskipamount}
\setenumerate[0]{leftmargin=*,itemsep=\the\smallskipamount}

\renewcommand{\to}{%
   \ifbool{@display}{\longrightarrow}{\rightarrow}%
   }
\let\shortmapsto\mapsto
\renewcommand{\mapsto}{%
   \ifbool{@display}{\longmapsto}{\shortmapsto}%
   }
\newlength{\olen}
\newlength{\ulen}
\newlength{\xlen}
\newcommand{\xra}[2][]{%
   \ifbool{@display}%
      {\settowidth{\olen}{$\overset{#2}{\longrightarrow}$}%
       \settowidth{\ulen}{$\underset{#1}{\longrightarrow}$}%
       \settowidth{\xlen}{$\xrightarrow[#1]{#2}$}%
       \ifdimgreater{\olen}{\xlen}%
          {\underset{#1}{\overset{#2}{\longrightarrow}}}%
          {\ifdimgreater{\ulen}{\xlen}%
             {\underset{#1}{\overset{#2}{\longrightarrow}}}
             {\xrightarrow[#1]{#2}}}}%
      {\xrightarrow[#1]{#2}}
   }
\makeatother
\newcommand{\xyra}[2][]{%
   \settowidth{\xlen}{$\xrightarrow[#1]{#2}$}%
   \ifbool{@display}%
      {\settowidth{\olen}{$\overset{#2}{\longrightarrow}$}%
       \settowidth{\ulen}{$\underset{#1}{\longrightarrow}$}%
       \ifdimgreater{\olen}{\xlen}%
          {\mathrel{\xymatrix@M=.12ex@C=3.2ex{\ar[r]^-{#2}_-{#1} &}}}%
          {\ifdimgreater{\ulen}{\xlen}%
             {\mathrel{\xymatrix@M=.12ex@C=3.2ex{\ar[r]^-{#2}_-{#1} &}}}
             {\mathrel{\xymatrix@M=.12ex@C=\the\xlen{\ar[r]^-{#2}_-{#1} &}}}}}%
      {\mathrel{\xymatrix@M=.12ex@C=\the\xlen{\ar[r]^-{#2}_-{#1} &}}}%
   }
\makeatletter
\newcommand{\xla}[2][]{%
   \ifbool{@display}%
      {\settowidth{\olen}{$\overset{#2}{\longleftarrow}$}%
       \settowidth{\ulen}{$\underset{#1}{\longleftarrow}$}%
       \settowidth{\xlen}{$\xleftarrow[#1]{#2}$}%
       \ifdimgreater{\olen}{\xlen}%
          {\underset{#1}{\overset{#2}{\longleftarrow}}}%
          {\ifdimgreater{\ulen}{\xlen}%
             {\underset{#1}{\overset{#2}{\longleftarrow}}}
             {\xleftarrow[#1]{#2}}}}%
      {\xleftarrow[#1]{#2}}
   }
\newcommand{\isoarrow}{%
   \ifbool{@display}{\overset{\sim}{\longrightarrow}}{\xrightarrow\sim}%
   }

\begin{document}

\title[]{Machine learning assisted exploration for affine Deligne-Lusztig varieties}

\author[Bin Dong]{Bin Dong}
\address{Beijing International Center for Mathematical Research, Peking University, Beijing, China; Center for Machine Learning Research, Peking University, Beijing, China}
\email{dongbin@math.pku.edu.cn}

\author[Xuhua He]{Xuhua He}
\address{Department of Mathematics and New Cornerstone Science Laboratory, The University of Hong Kong, Pokfulam, Hong Kong, Hong Kong SAR, China}
\email{xuhuahe@hku.hk}

\author[Pengfei Jin]{Pengfei Jin}
\address{Beijing International Center for Mathematical Research, Peking University, Beijing, China}
\email{jinpf@pku.edu.cn}

\author[Felix Schremmer]{Felix Schremmer}
\address{Department of Mathematics and New Cornerstone Science Laboratory, The University of Hong Kong, Pokfulam, Hong Kong, Hong Kong SAR, China}
\email{schremmer@hku.hk}

\author[Qingchao Yu]{Qingchao Yu}
\address{Beijing International Center for Mathematical Research, Beijing University, Haidian, Beijing, China}
\email{yuqingchao@bicmr.pku.edu.cn}

\thanks{}

\keywords{}
\subjclass[2010]{}

\date{\today}

\begin{abstract}

This paper presents a novel, interdisciplinary study that leverages a Machine Learning (ML) assisted framework to explore the geometry of affine Deligne-Lusztig varieties (ADLV). The primary objective is to investigate the nonemptiness pattern, dimension and enumeration of irreducible components of ADLV. Our proposed framework demonstrates a recursive pipeline of data generation, model training, pattern analysis, and human examination, presenting an intricate interplay between ML and pure mathematical research. Notably, our data-generation process is nuanced, emphasizing the selection of meaningful subsets and appropriate feature sets. We demonstrate that this framework has a potential to accelerate pure mathematical research, leading to the discovery of new conjectures and promising research directions that could otherwise take significant time to uncover. We rediscover the virtual dimension formula and provide a full mathematical proof of a newly identified problem concerning a certain lower bound of dimension. Furthermore, we extend an open invitation to the readers by providing the source code for computing ADLV and the ML models, promoting further explorations. This paper concludes by sharing valuable experiences and highlighting lessons learned from this collaboration.

\end{abstract}

\maketitle

\tableofcontents

\section{Introduction}

\subsection{A Brief Overview of Affine Deligne-Lusztig Varieties}\label{sec:adlvReview}

The concept of Affine Deligne-Lusztig Varieties (ADLV) was first introduced by Rapoport \cite{Rapoport2002}. These varieties serve as a group-theoretic model for the reduction of Shimura varieties and shtukas with parahoric level structure and play a vital role in arithmetic geometry and the Langlands program. Key problems associated with ADLV include:

\begin{itemize}
\item Non-emptiness pattern;
\item Dimension;
\item Enumeration of irreducible components.
\end{itemize}

Over the past two decades, the study of ADLV has been a vibrant research topic. Significant progress has been made in understanding fundamental problems, and important applications to number theory and the Langlands program have been discovered. The non-emptiness pattern and dimension of ADLV in the affine Grassmannian are now fully understood, and in most cases, they are also known for ADLV in the affine flag variety. The enumeration of irreducible components has been solved in the affine Grassmannian case. For more in-depth information, readers are referred to the survey article \cite{He2018}.

Despite these advancements, a comprehensive solution to the problems presented by ADLV remains a challenge, primarily due to the difficulty in finding explicit patterns.

In this paper, our focus is on ADLV $X_w(b)$ in the affine flag variety (associated with the Iwahori level structure). Information for other parahoric level structures can be obtained from the Iwahori level structure via the natural projection map. The ADLV $X_w(b)$ depends on two parameters: the element $w$ in the Iwahori-Weyl group $\tW$ of a loop group $\breve G$, and the Frobenius-twisted conjugacy class $[b]$ of $\breve G$. We consider the map from the pair $(w, b)$ to the dimension and enumeration of the top-dimensional irreducible components of $X_w(b)$ (refer to \S \ref{sec:pre-ADLV} for a precise statement of the problem). The image of a given pair $(w, b)$ can be calculated using the intricate inductive algorithm established in \cite{He2014_virtdim}. Yet, the goal is to derive more explicit formulas/information for this map. Such explicit formulas and information are particularly intriguing because of their broad applications to arithmetic geometry and number theory.

\subsection{Machine learning assisting pure mathematics research } 

In recent years, machine learning (ML), particularly deep learning, has exerted a profound impact on diverse scientific and engineering disciplines, bringing about substantial changes in the way we conduct research. The success of ML primarily hinges on the development of the designing and training of deep neural networks, which approximate complex, high-dimensional mappings with desirable accuracy and rapid evaluations. As a result, ML has become a leading force in different areas of artificial intelligence (AI), such as natural language processing (large language models such as chatGPT \cite{zhao2023survey}), computer vision (e.g., NerF \cite{mildenhall2021nerf}, Diffusion Models \cite{ho2020denoising}) and games (e.g., mastering the game of Go \cite{silver2016mastering}, Poker \cite{brown2019superhuman}, Starcraft II \cite{vinyals2019grandmaster}). Moreover, the excellent approximation capabilities of deep neural networks have helped discover new patterns or principles within large, multi-dimensional data sets. This has significantly expanded ML's role in natural science, contributing to the rise of a new field known as ``AI for Science." Successful examples include the work of AlphaFold \cite{jumper2021highly}, molecular dynamics simulations \cite{chmiela2017machine,zhang2018deep,jia2020pushing}, chemical discovery \cite{stokes2020deep,von2020retrospective,miethke2021towards}, system identification \cite{brunton2016discovering,long2018pde,long2019pde,raissi2019physics}, controllable nuclear fusion \cite{kates2019predicting,degrave2022magnetic}, etc.

More recently, Geordie Williamson and DeepMind used ML to assist in research-level explorations of pure mathematics \cite{davies2021advancing}. They present a ML-based framework that augments mathematicians' intuition, aiding in the discovery and understanding of complex mathematical relationships. This approach identifies potential correlations between two mathematical entities by deriving a function approximating the relationship and helping mathematicians analyze it.

The framework validates possible patterns in mathematical objects using supervised learning, and helps understand these patterns using attribution techniques. In the supervised learning stage, a hypothesis about a connection between two entities is proposed, a dataset is generated, and a function is trained to predict one entity from the other. The role of ML here is in learning a wide variety of potential non-linear functions given sufficient data. Attribution techniques are then used to understand the trained function and propose a potential relationship. One such technique, gradient saliency, calculates the derivative of function outputs with respect to inputs, helping to identify the most relevant problem aspects. This process may be iterative until a feasible conjecture is found.

In essence, this ML-guided framework enables a quick verification of the potential worthiness of an intuition about a relationship and, if validated, suggests how they may be related. This framework has already proven its utility by \cite{davies2021advancing} in achieving significant results, such as uncovering the first relationships between algebraic and geometric invariants in knot theory and conjecturing a resolution to the combinatorial invariance conjecture for symmetric groups.

\subsection{Our objective}

In this study, our objective is to develop a ML-assisted framework to guide the study of fundamental problems related to ADLV, specifically the nonemptiness pattern, the dimension and enumeration of irreducible components. As illustrated in Figures \ref{fig:pipeline} and \ref{fig:Program}, our framework showcases a recursive pipeline of data generation, model training, pattern analysis, and human examination. Despite similarities with the framework in the aforementioned study \cite{davies2021advancing}, several crucial differences exist. Our data-generation process is more intricate, particularly regarding the selection of a meaningful subset of pairs $(w,b)$ and an appropriate set of features. Additionally, after fitting a functional relationship between the feature set and the property of interest (e.g.\ the nonemptiness of $X_w(b)$), the patterns revealed by salience analysis may be more challenging for mathematicians to interpret due to the problem's complexity. Nevertheless, we found that this interaction between ML and human mathematicians significantly accelerates pure mathematical research, enabling us to identify new conjectures and promising research directions that could otherwise take years for mathematicians to discover by themselves. 

We provide the source code for computing geometric invariants of ADLV and machine learning models to invite interested readers to delve into this problem, reproduce our experiments, and refine our approach by studying different datasets and feature sets.
In Section \ref{sec:SearchingDimFormula}, we seek a linear approximation for the dimension of ADLV, leading us to rediscover the virtual dimension formula. Originally, the development of the virtual dimension formula required several years of intense research by many mathematicians and constitutes a major milestone in the field.
In Section \ref{sec:importantFeatures}, we conduct sensitivity analysis for the problems introduced in Section~\ref{sec:adlvReview}. We identify several important features, affirming some of the latest research results in the field and suggesting potential next steps.
Motivated by the experiments in Sections \ref{sec:SearchingDimFormula} and \ref{sec:importantFeatures}, we discover a new problem concerning a certain lower bound of the dimension, which has not been studied in the literature before. In Section \ref{sec:lowerBound}, we provide a full mathematical proof of this lower bound.
We conclude this paper by suggesting future work in Section \ref{sec:conclusion}, sharing some experiences, and highlighting lessons learned from the collaboration.

\noindent {\bf Acknowledgement: } BD and PJ are partially supported by NSFC 12090022. XH is partially supported by the New Cornerstone Science Foundation through the New Cornerstone Investigator Program and the Xplorer Prize, and by Hong Kong RGC grant 14300220. All authors are fully supported by their enthusiasm towards the emerging field of AI in Mathematics.

\section{Preliminaries}

\subsection{Definition and properties of affine Deligne-Lusztig varieties}\label{sec:pre-ADLV}
In this subsection, we provide a brief overview of affine Deligne-Lusztig varieties. Unless otherwise stated (i.e.\ sections \S \ref{sec:Program1} and \ref{sec:lowerBound}), we focus on the case of the special linear group $\SL_n$, which is also referred to as the type $A_{n-1}$. By focusing exclusively on this case, we may reduce the technical details in this exposition. However, the more important reason for this specialization is that it allows us to perform computer experiments within a reasonably narrow scope, where all the beauties and pathologies of the general case are still present.

Let $q$ be a prime power and $\BF_q$ be the finite field with $q$ elements. We define $F = \BF_q\doubleparenthesis t$ to be the field of formal Laurent series over $\BF_q$. This means that elements in $a\in F$ are formal power series
\begin{align*}
a = \sum_{i\in\mathbb Z} a_i t^i
\end{align*}
with coefficients $a_i\in \BF_q$ such that $a_i=0$ for almost all $i<0$. There is no notion of convergence involved, but the definition of addition and multiplication in $F$ mimics the behavior of absolutely convergent power series over real or complex numbers.

Pick once and for all an algebraic closure $\overline{\BF_q}$ and define $\breve F = \overline{\BF_q}\doubleparenthesis t$ to be the field of formal Laurent series over $\overline{\BF_q}$. The Galois group of the field extension $\breve F/F$ is generated by the Frobenius $\sigma$, which can be evaluated for elements in $\breve F$ as
\begin{align*}
\sigma\Bigl(\sum\nolimits_{i\in\mathbb Z} a_i t^i\Bigr) = \sum\nolimits_{i\in\mathbb Z}a^{q}_i t^i\in \breve F.
\end{align*}
Finally, we write $\mathcal O_{\breve F} = \overline{\BF_q}\doublebracket t$ for the ring of all formal power series, i.e.\ elements $a\in \breve F$ with $a_i=0$ for all $i<0$.

Throughout this paper till Section \ref{sec:lowerBound}, we will focus on the algebraic group scheme $\SL_n$. We get an induced map $\sigma : \SL_n(\breve F)\rightarrow\SL_n(\breve F)$, given by applying the above Frobenius $\sigma:\breve F\rightarrow \breve F$ to the entries of each $n\times n$-matrix in $\SL_n(\breve F)$.


Two elements $b,c\in \SL_n(\breve F)$ are called \emph{$\sigma$-conjugate} if there exists some $g\in \SL_n(\breve F)$ with
\begin{align*}
b = g^{-1} c \sigma(g).
\end{align*}
One checks that this is an equivalence relation, similar to ordinary conjugacy. We denote the $\sigma$-conjugacy class of $b\in \SL_n(\breve F)$ by $[b]$, and the set of $\sigma$-conjugacy classes of $\SL_n(\breve F)$ by $B(\SL_n)$. The $\sigma$-conjugacy class of $b\in\SL_n(\breve F)$ is uniquely determined by an invariant called the \emph{Newton point} of $b$, denoted $\nu_b\in \mathbb Q^n$ \cite{Kottwitz1985}.

Call a vector $(\nu_1,\dotsc,\nu_n)\in \mathbb Q^n$ \emph{dominant} if $\nu_1\ge\cdots\ge \nu_n$. Then the Newton point of each $b\in \SL_n(\breve F)$ is such a dominant vector. We have an action of the symmetric group $S_n$ on $\mathbb Q^n$ by permutation of coordinates. One checks that each orbit under this action contains precisely one dominant vector. If $b\in \SL_n(\breve F)$ is a diagonal matrix of the form $b = \diag(\pm t^{b_1},\dotsc,\pm t^{b_n})$ with $b_1,\dotsc,b_n\in\mathbb Z$, then the Newton point of $b$ is the unique dominant element in the $S_n$-orbit of $(b_1,\dotsc,b_n)\in\mathbb Z^n\subseteq \mathbb Q^n$.

One calls $W_0 := S_n$ the \emph{(finite) Weyl group} of $\SL_n$. The \emph{affine Weyl group} is given by the semidirect product
\begin{align*}
W_a = \tilde W := S_n\ltimes \{(\mu_1,\dotsc,\mu_n)\in\mathbb Z^n\mid \mu_1+\cdots + \mu_n=0\}.
\end{align*}
We write elements $w\in W_a$ also as $w = t^\lambda z$, where $z\in S_n$ and $\lambda =(\lambda_1,\dotsc,\lambda_n)\in\mathbb Z^n$ with $\lambda_1+\cdots+\lambda_n=0$.  The symbol $t$ is a formal variable reminding us of the uniformizer $t\in F$.

We choose for each permutation $x \in S_n$ a representative $\dot x \in \SL_n(\breve F))$ so that
\begin{align*}
\dot x_{i,j} = \begin{cases}\pm 1,&i=x(j),\\0,&i \neq x(j).\end{cases}
\end{align*}
So e.g.\ $\dot x$ may be chosen to be the permutation matrix of $x$ if $x$ is an even permutation, and the permutation matrix of $x$ with one sign flipped if $x$ is an odd permutation. For $w = t^\lambda z\in W_a$, we write $\dot w\in \SL_n(\breve F)$ for the element $\dot w =  \diag(t^{\lambda_1},\dotsc,t^{\lambda_n}) \dot z$.

Let $\breve I\subset \SL_n(\breve F)$ be the \emph{Iwahori} subgroup
\begin{align*}
\breve I =& \SL_n(\breve F)\cap \Bigl\{
\begin{pmatrix}\mathcal O_{\breve F}&\mathcal O_{\breve F}&\hdots&\mathcal \mathcal O_{\breve F}\\
t\mathcal O_{\breve F}&\mathcal O_{\breve F}&\ddots&\vdots\\
\vdots&&\ddots&\mathcal O_{\breve F}\\
t\mathcal O_{\breve F}&\hdots&t\mathcal O_{\breve F}&\mathcal O_{\breve F}\end{pmatrix}
\Bigr\}
\\=&\Bigl\{a\in \SL_n(\breve F)\mid \begin{array}{l}
\forall i> j:~a_{i,j}\in t\mathcal O_{\breve F}\\
\forall i\le j:~a_{i,j}\in \mathcal O_{\breve F}\end{array}
\Bigr\}
\\=&\{a\in \SL_n(\mathcal O_{\breve F})\mid a\text{ is upper triangular modulo }t\}.
\end{align*}
Then each element $g\in \SL_n(\breve F)$ has the form $g = i_1 \dot w i_2$ for a uniquely determined element $w\in W_a$ and non-unique $i_1, i_2\in \breve I$. Such a decomposition can be computed e.g.\ by using an adaption of the Gauss algorithm.

We have seen two decompositions of the set $\SL_n(\breve F)$, namely one into $\sigma$-conjugacy classes $B(\SL_n)$ and another one into Iwahori double cosets $\breve I\dot w \breve I$ for $w\in W_a$. We also write $\breve I w\breve I := \breve I\dot w \breve I$, since the double coset is independent of the choice of $\dot w$.

The right coset space $Fl=\SL_n(\breve F)/\breve I$ is called the \emph{affine flag variety}. It is an ind-scheme over $\overline{\BF_q}$, which behaves similarly to finite-dimensional varieties over that field. For $w\in W_a$ and $b\in \SL_n(\breve F)$, we define the \emph{affine Deligne-Lusztig variety} to be the subvariety
\begin{align*}
X_w(b) := \{g \breve I \in Fl\mid g^{-1} b \sigma(g)\in \breve Iw\breve I\}.
\end{align*}
Again, $X_w(b)$ is not actually a variety, but rather an ind-scheme over $\overline{\BF_q}$. However, each irreducible component of $X_w(b)$ is an actual finite-dimensional variety over $\overline{\BF_q}$. If $b$ is $\sigma$-conjugate to $c\in\SL_n(\breve F)$, say $c = h^{-1} b \sigma(h)$ for $h\in \SL_n(\breve F)$, then
\begin{align*}
X_w(b)\rightarrow X_w(c),\quad g\breve I\mapsto hg\breve I
\end{align*}
is an isomorphism. Thus we may associate the isomorphism type of $X_w(b)$ to the pair $(w,[b])\in W_a\times B(\SL_n)$. This is our main object of interest.

We see that the $\sigma$-centralizer of $b$, denoted
\begin{align*}
\bJ_b(F) = \{g\in \SL_n(\breve F)\mid g^{-1}b\sigma(g) = b\},
\end{align*}
acts on $X_w(b)$ by left multiplication. Up to that action, there are only finitely many irreducible components in $X_w(b)$. Since each such irreducible component is a finite-dimensional variety over $\overline{\BF_q}$, the entire affine Deligne-Lusztig variety $X_w(b)$ is always finite dimensional.

The main questions regarding the geometry of affine Deligne-Lusztig varieties are the following three:
\begin{itemize}
    \item Nonemptiness pattern: given $(w, [b])$, determine whether $X_w(b) \neq \emptyset$;
    \item Dimension: given $(w, [b])$ such that $X_w(b) \neq \emptyset$, calculate the dimension of $X_w(b)$;
    \item Irreducible components: given $(w, [b])$ such that $X_w(b) \neq \emptyset$, calculate the number of $\bJ_b(F)$-orbits of top dimensional irreducible components, i.e.\ the cardinality of $\bJ_b(F) \backslash \Sigma^{\mathrm{top}}(X_w(b))$. Here $\Sigma^{\mathrm{top}}(X)$ denotes the set of top dimensional irreducible components of $X$. 
\end{itemize}

\subsection{Important invariants}\label{sec:invariants}

To address these three questions, one explores the relationship between the affine Weyl group $W_a$ and the set $B(\SL_n)$ parametrized by the Newton points. This relationship is essentially combinatorial in nature, which allows us to compute the answers to the three above questions for any given pair $(w,\nu_b)$. We summarize some key combinatorial invariants that have proven valuable in previous works.

The group $W_0 = S_n$ is known to be a \emph{Coxeter group} with respect to the generators $\BS = \{s_1,\dotsc,s_{n-1}\}$. Here, $s_i$ is the simple reflection interchanging $i$ with $i+1$ and leaving everything else fixed. Each element $x\in W_0$ is a product of the $s_i$, and the shortest length of such an expression is called the \emph{length} of $x$. There is a unique element of maximal length in $W_0$, denoted $w_0$. It is the permutation $w_0(i) = n+1-i$ for $i\in\{1,\dotsc,n\}$, and its length is $\ell(w_0) = n(n-1)/2$.

There is a different way to compute the length of an element $x\in W_0$. Denote the set of \emph{roots} of $\SL_n$ as
\begin{align*}
\Phi = \{e_i - e_j\in\mathbb Q^n\mid i,j\in\{1,\dotsc,n\}\text{ and }i\neq j\}.
\end{align*}
The root $\alpha_{i,j}:=e_i-e_j$ is called \emph{positive} if $i<j$, and negative otherwise. We write $\delta(\alpha) = 1$ if $\alpha$ is a negative root and $0$ if $\alpha$ is positive. Observe that the $S_n$-action on $\mathbb Q^n$ preserves the set of roots. Then the length of $x\in W_0$ equals the number of positive roots $\alpha_{i,j}$ such that $x\alpha_{i,j}$ is a negative root. As formula,
\begin{align*}
\ell(x) = \sum_{i<j}\delta(x\alpha_{i,j}).
\end{align*}

The group $W_a = S_n\ltimes \mathbb Z^n$ is also known to be a Coxeter group with respect to $\BS_a$. The set of simple affine reflections $\BS_a$ is given by $\{s_0\}\cup \mathbb S$, where
\begin{align*}
s_0 = (1~n)t^{(-1,0,\dotsc,0,1)}\in W_a.
\end{align*}
One defines the length of an element $w\in W_a$ as above, given by the smallest representation using these simple affine reflections. There is an alternative way to compute the length of $w = t^{\lambda} z\in W_a$: We saw above that there is some $y\in W_0$ with $y z^{-1}\lambda\in \mathbb Z^n$ being dominant. Among all those elements in $y\in W_0$, there is a unique one with $\ell(y)$ being minimal. For this specific $y\in W_0$, we write $x := zy^{-1}\in W_0$ and $\mu := y z^{-1}\lambda\in\mathbb Z^n$, so that $w = xt^{\mu} y$. Then
\begin{align*}
\ell(xt^\mu y) = \langle \mu,2\rho\rangle + \ell(x)-\ell(y).
\end{align*}
Here, $\langle\cdot,\cdot\rangle$ is the standard Euclidean inner product on $\mathbb Q^n$, and $2\rho = (n-1,n-3,\dotsc,3-n,1-n)\in \mathbb Q^n$ is the sum of positive roots. Whenever we write $w$ in the form $xt^\mu y$, we always assume that $x,\mu,y$ have been chosen as above.

For $c\in\mathbb Z_{\ge 0}$, we call the element $w = xt^\mu y$ to be \emph{$c$-regular} if $\langle \mu,\alpha_{i,j}\rangle\ge c$ for all positive roots $\alpha_{i,j}$. The decomposition of $w$ into $x$, $t^\mu$, and $y$ has the most desirable properties whenever $w$ is $2$-regular, but the above length formula is always true even when such a regularity condition is not satisfied.

To each $w\in W_a$, one may associate the $\sigma$-conjugacy class $[\dot w]\in B(\SL_n)$. Its Newton point can be computed as follows: If $w$ has the form $w = t^\lambda$ for some $\lambda\in X_\ast(T)$, then $\dot w = \diag(t^{\lambda_1},\dotsc,t^{\lambda_n})$ and we saw above that the Newton point of $\dot w$ is the unique dominant element in the $S_n$-orbit of $\lambda$. For general $w\in W_a$, one may find an integer $m\ge 1$ such that $w^m$ is of the above form, and then the Newton point of $\dot w$ is given by $\nu_w = \nu_{w^m}/m\in\mathbb Q^n$.

It turns out that each $\sigma$-conjugacy class $[b]\in B(\SL_n)$ contains the representative $\dot w\in \SL_n(\breve F)$ of some $w\in W_a$, cf.\ \cite[Theorem~3.7]{He2014_virtdim}. Hence the above method yields all Newton points of all $\sigma$-conjugacy classes. If $[b]\in B(\SL_n)$ has Newton point $\nu_b = (\nu_1,\dotsc,\nu_n)\in\mathbb Q^n$, we define \emph{the best integral approximation} $\lfloor \nu_b\rfloor\in\mathbb Z^n$ to be the vector $(\mu_1,\dotsc,\mu_n)\in\mathbb Z^n$ such that for all $i$, we have
\begin{align*}
\lfloor \nu_1+\cdots + \nu_i\rfloor = \mu_1+\cdots + \mu_i\in\mathbb Z.
\end{align*}
Equivalently, $\lfloor \nu_b\rfloor$ is the unique vector in $\mathbb Z^n$ that can be written in the form
\begin{align*}
\lfloor \nu_b\rfloor = \nu_b - c_1 \alpha_{1,2}-\cdots - c_{n-1}\alpha_{n-1,n}
\end{align*}
such that $0\leq c_i<1$ for $i=1,\dotsc,n-1$. 
The \emph{defect} of $[b]\in B(\SL_n)$ can then be defined as $\de(b) = \langle \nu_b-\lfloor\nu_b\rfloor,2\rho\rangle$. Alternatively, it can be computed as the number of non-integral coordinates of $\nu$, i.e.\ the number of indices $i\in\{1,\dotsc,n\}$ such that $\nu_i\in\mathbb Q\setminus\mathbb Z$. 

\subsection{Computing the geometry of affine Deligne-Lusztig varieties}\label{sec:ADLVAlgo}

In this section, we present a combinatorial algorithm that efficiently computes the answers to the above three main questions for a given pair $(w,\nu_b)$, where $w\in W_a$ and $\nu_b\in\mathbb Q^n$. Although the subsequent sections of this article do not depend on the specific algorithm used, we want to provide a detailed explanation of its workings.

It is worth noting that the algorithm is somewhat complex and non-deterministic, which accounts for why the three main questions are still considered open. While the algorithm yields a computational solution to each of the three problems, one may still seek a more straightforward and satisfactory characterization. Nonetheless, the algorithm allows for effective computation of the desired results, which is crucial for our practical applications.

Let $w\in W_a$. We explain an algorithm that computes, in finite time, the set \begin{align*}\{\nu_b\mid X_w(b)\neq\emptyset\}\subseteq \mathbb Q^n\end{align*} (which hence must be finite). For each occurring Newton point, the corresponding affine Deligne--Lusztig variety is uniquely determined up to isomorphism, and our algorithm computes the dimension of this variety and the number of $\bJ_b(F)$-orbits of its top dimensional irreducible components.

For a simple affine reflection $s\in \BS_a$, we say that $sws\in W_a$ is a \emph{cyclic shift} of $w$ if $\ell(sws)\le \ell(w)$. Under this condition, $\ell(sws)$ can either be equal to $\ell(w)$ or $\ell(w)-2$.

In the first case, i.e.\ $\ell(sws) = \ell(w)$, the affine Deligne--Lusztig varieties $X_w(b)$ and $X_{sws}(b)$ are always isomorphic for all $[b]\in B(\SL_n)$. So in order to compute the above data for $w$, we may freely pass between $w$ and $sws$.

In the second case, $\ell(sws) = \ell(w)-2$, each affine Deligne--Lusztig variety $X_w(b)$ splits into two parts, so $X_w(b) = U\sqcup V$ is the disjoint union of two subsets, with $U$ being open in $X_w(b)$ and $V$ closed. Then there are surjective maps with irreducible one-dimensional fibres
\begin{align*}
U\twoheadrightarrow X_{ws}(b),\qquad V\twoheadrightarrow X_{sws}(b).
\end{align*}
Hence $U\neq \emptyset$ if and only if $X_{ws}(b)\neq \emptyset$. In this case, $\dim U = \dim X_{ws}(b)+1$. The set $U$ is $\bJ_b(F)$-invariant, and the number of $\bJ_b(F)$-orbits of top dimensional irreducible components agrees for $U$ and $X_{ws}(b)$. The same story happens for $V$ and $X_{sws}(b)$. Once we know this geometric information, the corresponding data for $X_w(b) = U\sqcup V$ is easily computed: If $U=V=\emptyset$, then $X_w(b)=\emptyset$. If precisely one of the subsets $U$ or $V$ is empty and the other one is non-empty, then $X_w(b)$ agrees with the unique non-empty subset, and all geometric invariants are known. Finally, if $U\neq\emptyset\neq V$, we have
$\dim X_w(b) =\max(\dim U,\dim V)$ and
\begin{align*}
\# \bJ_b(F)\setminus \Sigma^{\text{top}}(X_w(b)) =&\begin{cases} \#\bJ_b(F)\setminus \Sigma^{\text{top}}(U),&\dim U>\dim V,\\\#\bJ_b(F)\setminus \Sigma^{\text{top}}(V),&\dim V>\dim U,\\
 \#\bJ_b(F)\setminus \Sigma^{\text{top}}(U)+ \#\bJ_b(F)\setminus \Sigma^{\text{top}}(V),&\dim U=\dim V.\end{cases}
\end{align*}
Moreover, if $\ell(sws)=\ell(w)-2$, then $\ell(ws)=\ell(w)-1$. So in this case, we have reduced the geometric questions for $w$ and arbitrary $[b]$ to the same questions of the two elements $sws, ws$ of smaller length.

The first part of the algorithm iteratively enumerates the \emph{cyclic shift class} of $w$, i.e.\ the set of all elements in $W_a$ reachable by iterated cyclic shifts $w\rightarrow s_1 w s_1\rightarrow s_2 s_1 w s_1 s_2\rightarrow\cdots$. Each element in this cyclic shift class has length $\le \ell(w)$, so that the cyclic shift class is a finite set.

We traverse this cyclic shift class, until we either exhaust the full set or we reach some element $w'$ in the cyclic shift class and some $s'\in \BS_a$ with $\ell(w) = \ell(w') = \ell(s' w' s')+2$. In the latter case, we know $X_w(b)\cong X_{w'}(b)$ and the geometric properties of $X_{w'}(b)$ can be reduced to recursively calling our algorithm for the two smaller elements $s' w'$ and $s' w' s'$. Since the length drops by at least one, such a recursive call can only happen a finite number of times.

The second part of the algorithm handles the case where the entire cyclic shift class is enumerated without reaching a length-reducing cyclic shift. In this case, we not only know that $w$ has a minimal length in this cyclic shift class but even that it must have a minimal length in its conjugacy class in $W_a$.

In this case, we can explicitly state that
\begin{align*}
\{\nu_b\mid X_w(b)\neq\emptyset\} = \{\nu_w\}.
\end{align*}
For $\nu_b = \nu_w$, we know $\dim X_w(b) = \ell(w) - \langle \nu_b,2\rho\rangle$ and that the number of $\bJ_b(F)$-orbits of top dimensional irreducible components is equal to $1$.

This algorithm is guaranteed to terminate in finite time with the correct result. We note that the algorithm itself is non-deterministic since there is no canonical way to traverse the cyclic shift class of an element $w\in W_a$. While the above method allows us to compute the dimension of $X_w(b)$ for arbitrary $w,[b]$, it is far from being a closed formula. In many cases, however, one may expect that such a closed formula can be found. 

The first part of this algorithm is due to G\"ortz-He \cite[Corollary~2.5.3]{Goertz2010}, the second part is due to He and Nie \cite[Theorem~A]{He2014_minlength}, \cite[Theorem~4.8]{He2014_virtdim}.

\subsection{Machine Learning assisted Formula Exploration}

As mentioned earlier, the algorithm used to compute the geometry of affine Deligne-Lusztig varieties is somewhat complex, non-deterministic, and implicit. However, practical research often requires an explicit expression or pattern. Machine learning, particularly deep neural network models, excels at  fitting and analyzing high-dimensional mappings. Thus, we are considering utilizing machine learning to aid us in exploring formulas.

\begin{figure}[htb]
\centering
\begin{tikzpicture}[x=0.75pt,y=0.75pt,yscale=-1,xscale=1]

\draw [line width=1.5]    (259,64) -- (310.67,64) ;
\draw [shift={(313.67,64)}, rotate = 180] [color={rgb, 255:red, 0; green, 0; blue, 0 }  ][line width=1.5]    (14.21,-4.28) .. controls (9.04,-1.82) and (4.3,-0.39) .. (0,0) .. controls (4.3,0.39) and (9.04,1.82) .. (14.21,4.28)   ;
\draw [line width=1.5]    (404,95) -- (404,123.1) ;
\draw [shift={(404,126.1)}, rotate = 270] [color={rgb, 255:red, 0; green, 0; blue, 0 }  ][line width=1.5]    (14.21,-4.28) .. controls (9.04,-1.82) and (4.3,-0.39) .. (0,0) .. controls (4.3,0.39) and (9.04,1.82) .. (14.21,4.28)   ;
\draw [line width=1.5]    (321,169) -- (259.65,169) ;
\draw [shift={(256.65,169)}, rotate = 360] [color={rgb, 255:red, 0; green, 0; blue, 0 }  ][line width=1.5]    (14.21,-4.28) .. controls (9.04,-1.82) and (4.3,-0.39) .. (0,0) .. controls (4.3,0.39) and (9.04,1.82) .. (14.21,4.28)   ;
\draw   (83.65,26.22) -- (235.65,26.22) -- (235.65,84.1) -- (83.65,84.1) -- cycle ;
\draw   (341.65,25.22) -- (493.65,25.22) -- (493.65,83.1) -- (341.65,83.1) -- cycle ;
\draw   (82.65,140.22) -- (234.65,140.22) -- (234.65,198.1) -- (82.65,198.1) -- cycle ;
\draw   (342.65,138.22) -- (494.65,138.22) -- (494.65,196.1) -- (342.65,196.1) -- cycle ;
\draw    (263.47,105.93) .. controls (316.2,51.21) and (327.36,173.7) .. (263.44,122.72) ;
\draw [shift={(262.47,121.93)}, rotate = 39.19] [color={rgb, 255:red, 0; green, 0; blue, 0 }  ][line width=0.75]    (10.93,-3.29) .. controls (6.95,-1.4) and (3.31,-0.3) .. (0,0) .. controls (3.31,0.3) and (6.95,1.4) .. (10.93,3.29)   ;
\draw [line width=1.5]  [dash pattern={on 5.63pt off 4.5pt}]  (161.47,127.93) -- (161.63,98) ;
\draw [shift={(161.65,95)}, rotate = 90.32] [color={rgb, 255:red, 0; green, 0; blue, 0 }  ][line width=1.5]    (14.21,-4.28) .. controls (9.04,-1.82) and (4.3,-0.39) .. (0,0) .. controls (4.3,0.39) and (9.04,1.82) .. (14.21,4.28)   ;

\draw (354.65,38) node [anchor=north west][inner sep=0.75pt]  [font=\normalsize] [align=left] {Generated Data \\ \ \ \ \ $\displaystyle \{X_{i} ,Y_{i}\}_{i=1}^{N}$};
\draw (358.67,158) node [anchor=north west][inner sep=0.75pt]  [font=\small] [align=left] {The Approximation\\ \ \ \ \ \ \ \ \ \ \ \ \ \ $\displaystyle \hat{f}_{\theta ^{*}}$};
\draw (115,156) node [anchor=north west][inner sep=0.75pt]  [font=\normalsize] [align=left] {Patterns of \\ \ \ \ \ \ \ \ \ $\displaystyle f$};
\draw (90.65,42.22) node [anchor=north west][inner sep=0.75pt]  [font=\normalsize] [align=left] {Concerned mapping \\$\displaystyle Y=f(X)$};
\draw (256,26) node [anchor=north west][inner sep=0.75pt]  [font=\normalsize] [align=left] {Generate\\Data};
\draw (412,97) node [anchor=north west][inner sep=0.75pt]  [font=\normalsize] [align=left] {Train a \\model $\displaystyle \hat{f}_{\theta}$};
\draw (263,179) node [anchor=north west][inner sep=0.75pt]  [font=\normalsize] [align=left] {Analysis\\of $\displaystyle \hat{f}_{\theta^{*}}$};
\draw (73,107) node [anchor=north west][inner sep=0.75pt]  [font=\normalsize] [align=left] {Modify $\displaystyle f(X)$};

\end{tikzpicture}
\caption{Pipeline of Machine Learning assisted Formula Exploration}\label{fig:pipeline}
\end{figure}
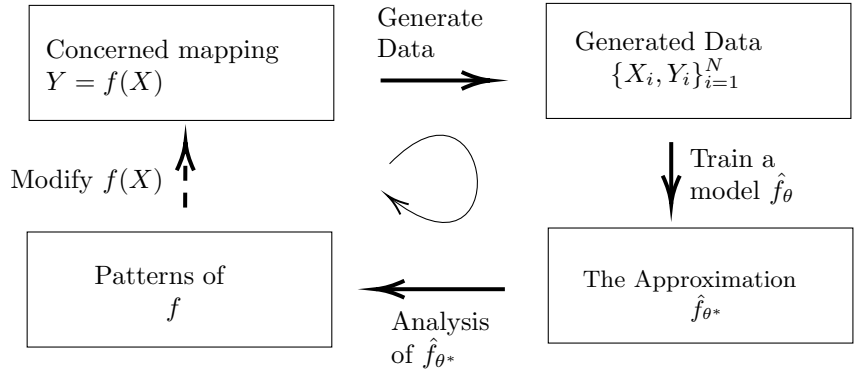

The machine learning-assisted formula exploration process is broken down into several steps, as illustrated in Figure \ref{fig:pipeline}. First, a suitable problem $Y=f(X)$ needs to be selected, where $f$ is a mapping from known features $X$ to the variable $Y$ of interest, and a dataset $\{X_i,Y_i\}_{i=1}^N$ is generated using a known algorithm. Second, a machine learning model $\hat{f}_\theta$ is chosen, and an optimal approximation $\hat{f}_{\theta^*}$ of $f$ is attained through optimization on the generated dataset. Third, hints about the patterns of $f$ are obtained by analyzing the explicit expression of $\hat{f}_{\theta^*}$. Finally, the obtained pattern can inspire us to introduce new features $X$ or consider functions $f$ with different domains. From this, we can modify $f(X)$, return to the first step, and continue this cycle. Each step will be elaborated on in detail in the following.

Regarding $f$, our primary concern is mapping form $(w, [b])$ to the geometry of affine Deligne-Lusztig varieties $X_w(b)$, including dimension, whether $X_w(b) \neq \emptyset$ and number of irreducible components. However, due to the inherent complexity of this mapping, it may be necessary to use the computed features of $w$ and $b$ as $X$, such as $\ell(w)$ and $\text{defect}(b)$.

The selection of a suitable machine learning model $\hat{f}_{\theta}$ requires careful consideration of the inherent complexity of $f$. Generally speaking, complex models have greater expressive power, but they may be more difficult to optimize and may require larger datasets. In this article, we focus on several commonly used models, including linear models \cite{aitken1936iv, nelder1972generalized}, Support Vector Machine (SVM) \cite{cortes1995support, burges1998tutorial}, and neural networks \cite{rosenblatt1958perceptron, rumelhart1986learning}. The specific formulation of the models will be described in detail in the subsequent section.


Given a dataset $\mathcal{D}=\{X_i,Y_i\}_{i=1}^N$, and the model $\hat{f}_{\theta}$, the next step is to identify the optimal value for $\theta$ that minimizes the distance between $\hat{f}_{\theta}$ and the target function $f$ (known as training) and evaluate the performance of $\hat{f}_{\theta}$ (known as testing). Typically, the dataset $\mathcal{D}$ is partitioned into two distinct subsets: the training set $\mathcal{D}_{tr}$ and the testing set $\mathcal{D}_{te}$. The former is used to obtain the optimal $\theta^*$, while the latter is used to evaluate $f_{\theta^*}$.

The training process seeks to minimize the following loss function:
$$\theta^*=\arg\min_\theta \sum_{{(X_i,Y_i)}\in\mathcal{D}{tr}}\mathcal{L}(Y_i,\hat{f}_{\theta}(X_i))+\mathcal{R}(\theta),$$
where $\mathcal{L}$ is a distance metric function, such as cross-entropy for classification problems or $L_2$ distance for regression problems. The regularization term for the parameters, $\mathcal{R}(\theta)$, depends on prior knowledge about the parameters and is typically expressed using the $L_1$ and $L_2$ norm to ensure simplicity or sparsity of the expression and prevent overfitting \cite{bach2011convex, loshchilov2017decoupled}.

After obtaining $\theta^*$, we typically calculate the loss function and accuracy on the testing set. If the model $\hat{f}_{\theta^*}$ exhibits relatively low loss and high accuracy on both the training and testing sets, we can consider it as a good approximation of the target function $f$. This indicates that the model has successfully generalized from the training set to unseen data and is likely to perform well on new data.

Once we obtain $\hat{f}_{\theta^*}$, we analyze the patterns of $f$ by examining $\hat{f}_{\theta^*}$. Firstly, the complexity of $f$ can be analyzed by observing the accuracy of $\hat{f}_{\theta^*}$ under different hyperparameters, such as the number of hidden neurons and layers in the neural network. This allows for a rough estimation of $f$'s complexity. Secondly, we can determine the sensitivity of $f$ to different features by taking the derivative of $\hat{f}_{\theta^*}$, enabling the determination of the significance of these features. Thirdly, if the form of $\hat{f}_{\theta^*}$ is relatively simple or becomes simple after sparse optimization, an approximate explicit expression of $f$ can be directly obtained. Finally, error analysis of $f_{\theta^*}$ also facilitates the understanding of $f$'s properties, such as differences in complexity in varying regions.

In cases where a suitable pattern cannot be obtained, this could indicate an improper selection of our mapping. For instance, if the number of features $X$ is insufficient or if $f$ is too complex, it may be challenging to find a simple $\hat{f}_\theta$ that can approximate the function and reveal its underlying patterns. In such instances, we need to modify $f(X)$ based on the insights gained from the previous round of exploration, such as by adding specific features or considering the properties of $f$ on specific domains, and continue with the next round of exploration. This process continues iteratively. 

\section{Fundamental Concepts of Machine Learning and Associated Caveats}

Machine learning is a field that employs computational models to learn patterns in data. This section will provide an overview of some basic machine learning models and discuss several crucial caveats in employing these models.

\subsection{Machine Learning Models}

\subsubsection{Linear Models} 
\label{sec:LinR}

Linear models are perhaps the simplest type of machine learning model, and they make a good starting point for the study of machine learning algorithms. They model the relationship between the input features and the output as a linear combination of the input features. 

Suppose we have \(p\) input features, a linear model is a hyperplane and is given by the equation:
\[\hat{Y}= \hat{f}_{\theta}(X) = \beta^\top X-b = \beta_{(1)} X_{(1)} + \beta_{(2)}X_{(2)} + ... + \beta_{(p)} X_{(p)}-b.\]
Here, \([X_{(1)}, X_{(2)}, ..., X_{(p)}]=X\in \mathbb{R}^p\) are the input features, \(\hat{Y}\) is the output, and \(\{[\beta_{(1)},...,\beta_{(p)}]=\beta\in \mathbb{R}^p, b\}=\theta\) are the parameters of the model. The parameters are typically learned from the data using a method called \textit{least squares} which minimizes minimizes the sum of the squared residuals, the differences between the observed data $Y_i$ and predicted output $\hat{Y}_i$:
\[\min_{\beta, b} \sum_{i=1}^{N} (Y_i - \beta^\top X_i +b )^2.\]
Despite their simplicity, linear models can be quite effective in practice, particularly when the data is actually linearly separable or close to it.

\subsubsection{Support Vector Machines (SVM)}
\label{sec:LinC}

Support Vector Machines (SVMs) are a set of supervised learning methods particularly well-suited for classification of complex but small or medium-sized datasets.

Fundamentally, SVMs aim to find a hyperplane that best separates the classes in the feature space by maximizing the margin. The margin is defined as the distance from the hyperplane to the nearest data point from either class. In the case of a linear SVM, the goal is to find the optimal hyperplane that maximizes this margin. The equation of this hyperplane is given by:
\[\beta^\top  X - b = 0\]

The optimization problem of finding the best hyperplane can be expressed as follows:
\[\min_{\beta, b} \frac{1}{2} ||\beta||^2,\quad Y_i (\beta^\top X_i - b) \geq 1, \quad \forall i,\]
where \(Y_i\) are the class labels and \(X_i\) are the data points.

For non-linearly separable data, SVMs utilize the ``kernel trick”, a method to map the input data into a higher-dimensional space where it can be linearly separable. Different kernel functions can be used depending on the nature of the data, such as polynomial kernels and radial basis function (RBF) kernels.
Despite their mathematical complexity, SVMs have a geometric interpretation and can be intuitively understood as trying to find the ``widest possible street" that separates the different classes.

\subsubsection{Neural Network (NN) Regression}
\label{sec:NnR}

Neural networks are flexible function approximators that use layers of neurons to model complex relationships. In a regression context, a neural network learns a mapping from inputs to a continuous output. Each neuron applies a series of transformations, first a linear transformation and then a
non-linear activation function, to its input.

Consider a fully connected network \cite{nair2010rectified} with a total of $n_L$ layers, each containing $n_H$ neurons, where the input is $X^{(0)}=X$ and the output is $\hat{Y}=\hat{f}(X)=X^{(n_L)}$. The neuron \(j\) in layer \(l\) can be represented as:
\[X_{(j)}^{(l)} = \rho\left(\sum_{i=1}^{n_H} w_{ji}^{(l)} X_{(i)}^{(l-1)} + b_{j}^{(l)}\right)\]
where \(w_{ji}^{(l)}\) and \(b_j^{(l)}\) are the weights and bias for neuron \(j\) in layer \(l\), \(X_{(i)}^{(l-1)}\) are the outputs of the neurons in the previous layer, and $\rho$ is the activation function. Common choices for $\rho$ include the sigmoid, hyperbolic tangent, and rectified linear unit (ReLU) functions. The network's weights and biases are learned by minimizing a loss function, often the mean squared error for regression tasks, using an algorithm such as gradient descent or one of its variants.

\subsubsection{Neural Network (NN) Classification}
\label{sec:NnC}

In a classification context, a neural network learns to classify inputs into discrete categories. The architecture is similar to that of a regression network, but the final layer typically has as many neurons as there are classes, and it applies a softmax function to produce a probability distribution $p$ over the classes.

The output of the \(j\)th neuron in the softmax layer, denoted as $p_{(j)}$, can be represented as:
\[p_{(j)} = \frac{e^{X_{(j)}^{(n_L)}}}{\sum_{k=1}^{K} e^{X_{(k)}^{(n_L)}}}\]
where \(X_{(j)}^{n_L}\) is the input to the \(j\)th neuron in the softmax layer, and \(K\) is the number of classes. The weights and biases are optimized by minimizing the cross-entropy loss, which measures the discrepancy between the predicted and true probability distributions across the classes.

\subsection{Caveats}

Machine learning methods have demonstrated tremendous efficacy across a variety of tasks. However, several potential caveats may affect their performance and interpretation:

\subsubsection{Choice of model:}\label{sec:caveatsChoiceOfModel} In a scenario where your primary goal is salience analysis, the choice of a machine learning model is heavily influenced by interpretability, ability to reveal feature importance, and the potential to infer functional forms between inputs and outputs. Here are some considerations to help guide your model selection:

\begin{enumerate}
\item {Interpretability:} When the aim is to understand what input features are important and their relationships with the output, interpretability becomes a crucial criterion for model selection. Models such as linear regression, logistic regression, and decision trees are traditionally more interpretable than, say, deep neural networks or support vector machines.

\item {Inherent Feature Importance:} Some models inherently provide feature importance measures. For example, the size of coefficients in linear models can indicate feature importance, and tree-based models provide feature importance based on the frequency of a feature being used to split the data.

\item {Flexibility v.s. Interpretability:} More flexible models like neural networks can model complex relationships, but they often lack interpretability. Conversely, simpler models like linear regression provide clearer insight into relationships between variables but may fail to capture complex, nonlinear relationships.

\item {Trade-off between Accuracy and Interpretability:} There is often a trade-off between model accuracy and interpretability. In salience analysis, we might be willing to sacrifice some accuracy for better interpretability, which should be factored into the model selection process.

\end{enumerate}
    
\subsubsection{Training and testing, overfitting:} One common pitfall in machine learning is overfitting, where the model learns the training data too well and performs poorly on unseen test data. In salience analysis scenarios, understanding the relationship between features and output is crucial. To ensure the model captures the actual underlying relationships and not the noise, controlling overfitting is essential. Here are some strategies that might help:

\begin{enumerate}
\item {Regularization:} This technique penalizes the complexity of the model, discouraging learning overly complex patterns that might be due to noise. The common forms of regularization include $\ell_1$- and $\ell_2$-regularization.

\item {Early Stopping (for neural networks):} During the training process, monitor the model's performance on a validation set. Stop training as soon as the performance on the validation set begins to degrade.

\item {Dropout (for neural networks):} Randomly "dropping out" units in a neural network during training can prevent complex co-adaptations on training data, which helps to avoid overfitting \cite{srivastava2014dropout}.

\item {Interpretable Models:} If the main goal is to understand the relationships between features and output, using simpler, interpretable models such as linear regression or decision trees could be beneficial. These models may be less prone to overfitting compared to complex models like deep neural networks.
\end{enumerate}

\subsubsection{Randomness:} When conducting salience analysis with machine learning models, it's crucial to understand and handle the randomness introduced by stochastic training processes. You want to ensure that the salience you identify is not due to the randomness in the training process, but truly indicative of the underlying relationships in your data. Below are a few strategies specifically for this context:

\begin{enumerate}
\item {Feature Importance Measures:} Many models provide ways to measure the importance of features, either directly (like coefficients in linear models), or indirectly (like gradients in neural networks). However, remember that these measures are impacted by the stochasticity of training. To mitigate this, consider averaging feature importance over multiple runs or training multiple models using different random seeds and averaging their feature importance measures.

\item {Ablation Studies:} One way to understand the importance of a feature is to see how much the model's performance drops when that feature is removed. By conducting this analysis over multiple runs (with different random seeds), you can obtain a measure of feature importance that is robust to the randomness of the training process.


\item {Controlled Training Processes:} Reduce the randomness in the training process through techniques such as decreasing the learning rate over time, using a larger batch size, or using a different optimization algorithm less sensitive to stochasticity, like RMSprop \cite{ruder2016overview} or Adam \cite{kingma2014adam}.
\end{enumerate}

\subsubsection{Data imbalance:} Machine learning models can perform poorly when there's a class imbalance in the training data. In such cases, the model might be biased towards the majority class. Handling imbalanced data is crucial, especially in a context where understanding the relationships between features and output is of primary importance. We list some strategies specifically geared towards such scenarios:

\begin{enumerate}

\item {Resampling Techniques:} One can alter the dataset itself to address the imbalance:
    \begin{itemize}
        \item \textit{Oversampling the Minority Class:} This involves creating or synthesizing new instances of the minority class until it reaches a similar number as the majority class. While it can balance the classes, it might lead to overfitting due to the replication of the minority instances.
        \item \textit{Undersampling the Majority Class:} This involves removing instances from the majority class until it reaches a similar number as the minority class. While it can be effective in balancing the dataset, it may cause loss of information by excluding potentially important instances from the majority class.
    \end{itemize}
    Both approaches aim to balance the distribution between the majority and minority classes but come with potential drawbacks that must be carefully considered during implementation.

\item {Cost-Sensitive Learning:} While this technique can always be used in cases of imbalanced data, the costs need to be chosen carefully in this context. A simple heuristic like setting the cost inversely proportional to class frequency might not be the best choice, as it might lead to the model focusing too much on rare classes that might not have enough instances to derive reliable feature-output relationships.

Another method is adjusting the temperature of the softmax function in the output layer of the model. The softmax function is often used in the final layer of a classification neural network to convert the outputs to probability values for each class. You can adjust the temperature based on the class frequencies, so that the model is made more sensitive to the minority class. However, keep in mind that it needs to be used carefully, as it can make the model more prone to overfitting to the minority class.
\end{enumerate}

\subsubsection{Dataset size: } Throughout this paper, we engage with datasets of different sizes depending on the complexity of the model under study - ranging from simple linear models to more elaborate two or three-layer neural networks. When the goal is to extract a concise and clear formula, we prefer simpler models such as linear models, which have fewer parameters, thereby requiring less data. On the other hand, for feature analysis, a more accurate representation of the unknown function relationship is desired, without overfitting, which calls for more sophisticated models and larger datasets. In practice, the actual size of a dataset is typically determined through an empirical approach that involves training a model on datasets of different sizes and comparing the results. Upon reaching a point where improvement is only marginal, we consider the dataset size to be sufficient. 

\section{Programs}

\begin{figure}[htb]
\centering
\begin{tikzpicture}[x=0.75pt,y=0.75pt,yscale=-1,xscale=1]

\draw   (88,41) -- (186.6,41) -- (186.6,102.6) -- (88,102.6) -- cycle ;
\draw   (281,42) -- (361.65,42) -- (361.65,96.23) -- (281,96.23) -- cycle ;
\draw   (462,45.23) -- (544.1,45.23) -- (544.1,88.37) -- (462,88.37) -- cycle ;
\draw   (72.43,151) -- (230.65,151) -- (230.65,221.3) -- (72.43,221.3) -- cycle ;
\draw   (461.1,140) -- (555.1,140) -- (555.1,222.3) -- (461.1,222.3) -- cycle ;
\draw    (20,72) -- (66.43,72) ;
\draw [shift={(68.43,72)}, rotate = 180] [color={rgb, 255:red, 0; green, 0; blue, 0 }  ][line width=0.75]    (10.93,-3.29) .. controls (6.95,-1.4) and (3.31,-0.3) .. (0,0) .. controls (3.31,0.3) and (6.95,1.4) .. (10.93,3.29)   ;
\draw    (209.43,71) -- (262.43,71) ;
\draw [shift={(264.43,71)}, rotate = 180] [color={rgb, 255:red, 0; green, 0; blue, 0 }  ][line width=0.75]    (10.93,-3.29) .. controls (6.95,-1.4) and (3.31,-0.3) .. (0,0) .. controls (3.31,0.3) and (6.95,1.4) .. (10.93,3.29)   ;
\draw    (373.65,72) -- (432.43,72) ;
\draw [shift={(434.43,72)}, rotate = 180] [color={rgb, 255:red, 0; green, 0; blue, 0 }  ][line width=0.75]    (10.93,-3.29) .. controls (6.95,-1.4) and (3.31,-0.3) .. (0,0) .. controls (3.31,0.3) and (6.95,1.4) .. (10.93,3.29)   ;
\draw    (135,113) -- (135.41,146.3) ;
\draw [shift={(135.43,148.3)}, rotate = 269.3] [color={rgb, 255:red, 0; green, 0; blue, 0 }  ][line width=0.75]    (10.93,-3.29) .. controls (6.95,-1.4) and (3.31,-0.3) .. (0,0) .. controls (3.31,0.3) and (6.95,1.4) .. (10.93,3.29)   ;
\draw    (502,97) -- (502.41,130.3) ;
\draw [shift={(502.43,132.3)}, rotate = 269.3] [color={rgb, 255:red, 0; green, 0; blue, 0 }  ][line width=0.75]    (10.93,-3.29) .. controls (6.95,-1.4) and (3.31,-0.3) .. (0,0) .. controls (3.31,0.3) and (6.95,1.4) .. (10.93,3.29)   ;
\draw    (245,193) -- (318.65,193) ;
\draw    (318.65,193) -- (318.65,162.23) -- (318.65,125.23) ;
\draw [shift={(318.65,123.23)}, rotate = 90] [color={rgb, 255:red, 0; green, 0; blue, 0 }  ][line width=0.75]    (10.93,-3.29) .. controls (6.95,-1.4) and (3.31,-0.3) .. (0,0) .. controls (3.31,0.3) and (6.95,1.4) .. (10.93,3.29)   ;
\draw  [dash pattern={on 2.53pt off 3.02pt}][line width=2.25]  (11,14) -- (271.65,14) -- (271.65,235.23) -- (11,235.23) -- cycle ;
\draw  [dash pattern={on 2.53pt off 3.02pt}][line width=2.25]  (368.1,12.57) -- (582.65,12.57) -- (582.65,235.23) -- (368.1,235.23) -- cycle ;

\draw (19,29) node [anchor=north west][inner sep=0.75pt]   [align=left] {Choose \\a Group};
\draw (97,66) node [anchor=north west][inner sep=0.75pt]   [align=left] {Type $\displaystyle A-G$};
\draw (204,31) node [anchor=north west][inner sep=0.75pt]   [align=left] {Generate\\Data};
\draw (284,65) node [anchor=north west][inner sep=0.75pt]   [align=left] {$\displaystyle \{X_{i} ,Y_{i}\}_{i=1}^{N}$};
\draw (376,32) node [anchor=north west][inner sep=0.75pt]   [align=left] {Training\\Model};
\draw (491,57) node [anchor=north west][inner sep=0.75pt]   [align=left] {$\displaystyle \hat{f}_{\theta ^{*}}$};
\draw (140,121) node [anchor=north west][inner sep=0.75pt]   [align=left] {Compute};
\draw (512,107) node [anchor=north west][inner sep=0.75pt]   [align=left] {Analysis};
\draw (80,163) node [anchor=north west][inner sep=0.75pt]   [align=left] {Dimension and \\Number of irreducible \\components of ADLV};
\draw (470.43,146) node [anchor=north west][inner sep=0.75pt]   [align=left] {Accuracy,\\Sensitivity,\\Visualization,\\.......};
\draw (282,172) node [anchor=north west][inner sep=0.75pt]   [align=left] {$\displaystyle Y_{i}$};
\draw (69,245) node [anchor=north west][inner sep=0.75pt]   [align=left] {Traditional Algorithms};
\draw (398,242.67) node [anchor=north west][inner sep=0.75pt]   [align=left] {Machine Learning Models};

\end{tikzpicture}

\caption{Workflow of the Program}\label{fig:Program}
\end{figure}
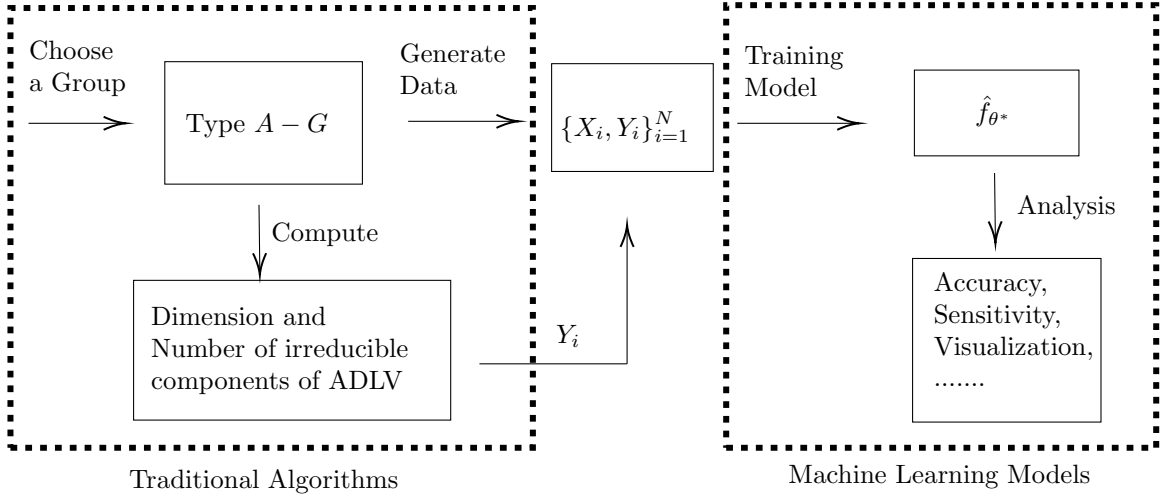

In this section, we introduce our program, which is composed of two primary modules. The first module computes variables of interest in affine Deligne-Lusztig varieties, as described in section \ref{sec:Program1}. The second module analyzes the data generated using machine learning techniques, as detailed in section \ref{sec:Program2}. Figure \ref{fig:Program} illustrates the overall workflow of the program.

Both parts of the program are publicly available online: \url{https://github.com/Jinpf314/ML4ADLV/}.

\subsection{Program for affine Deligne-Lusztig varieties}
\label{sec:Program1}
We give a short introduction to the program for the affine Deligne-Lusztig varieties based on Python. 

\textbf{Choosing a group.} The first input determines the type of algebraic group. The valid inputs are An, Bn, Cn, Dn, 2An, 2Dn.

To simplify notation, we focus on the group $G=SL_n$ (type $A_{n-1}$) throughout this section. However, we note that the program is compatible with all classical groups.

\textbf{Input for $w$.} The element $w \in W_a$ can be expressed in the following two ways:
\begin{itemize}
    \item The product of the translation part and the finite part; 
    \item The product of a sequence of simple reflections. 
\end{itemize}
Element of the form $w = t^{(a_1,\ldots,a_n)}s_{i_1} \cdots s_{i_r}$, where $(a_1,\ldots , a_n)\in X_*$ and all $i_j \in \{1,2,\ldots,n-1\}$, is written as $\text{affine}\_\text{Weyl}([a_1\ldots,a_n],[i_1,\dots,i_r])$.
Elements of the form $w =  s_{i_1} \cdots s_{i_r}$, where all $i_j \in \{0,1,2,\ldots,n-1\}$ is written as \text{exp}($[i_1\ldots,i_r])$. For example, in n=3 case (type A2). affine$\_$Weyl([1, 0, -1],[1, 2])=exp([0, 2]).

The simple reflections are s[0], s[1],$\ldots$, s[n-1]. The identity element of $W_a$ is Id.

\textbf{Input for $b$.} In this program, we input the Newton point $\nu_b$ of $b$ instead of $b$. Note that in the case of $\bG = \SL_n$, the conjugacy class $[b]$ is determined by its Newton point.






\begin{function*}[\textbf{dim}]Computing dimension of ADLV:
$$\text{dim} (w,\nu) $$\par
\textbf{Input: }$w\in\tW; \nu\in\mathbb{Q}^{n}$; \par
\textbf{Output: }$\dim X_w(b) $\par
\textbf{Description: }
If $X_w(b)=\emptyset$, the output is `empty'.

\end{function*}

\begin{function*}[\textbf{irr}]Computing irreducible components of ADLV:
$$\text{irr}(w,\nu)$$\par
\textbf{Input: }$w\in\tW; \nu\in\mathbb{Q}^{n}$; \par
\textbf{Output: }$ \sharp \bJ_b(F) \backslash \Sigma^{top}X_w(b)$\par
\textbf{Description: }
If $X_w(b)=\emptyset$, the output is 0.

\end{function*}

\begin{function*}[\textbf{$dim\_irr\_print$}]Listing all $b$ such that $X_w(b)\ne\emptyset$ and computing dimension and number of irreducible components:
$$\text{dim}\_\text{irr}\_\text{print}(w)$$\par
\textbf{Input: }$w\in\tW$ ; \par
\textbf{Output: }print the following:\par
Newton point = $\nu$, dim = $\text{dim}X_w(b)$, irr = $\sharp \bJ_b(F)\backslash \Sigma^{top}X_w(b)$

\textbf{Description:}
The function lists all $b$ such that $X_w(b)\ne\emptyset$.

\end{function*}

\begin{example*}
A2 case.\par
\textbf{Input: }\par 
w = affine$\_$Weyl([1,1,-2],[2,1]); dim$\_$irr$\_$print(w) \par
\textbf{Ouput: }\par
Newton point = [1/2, 1/2, -1], dim = 1, irr = 1\par
Newton point = [0, 0, 0], dim = 3, irr = 1\par
\textbf{Input: }\par
print(dim(w,[0,0,0]), irr(w,[1/2,1/2,-1]), dim(w,[1,0,-1]), irr(w,[2,0,-2]))\par
\textbf{Ouput: }\par
3 1 empty 0
\end{example*}

\subsection{Program for Machine Learning}
\label{sec:Program2}
The program implemented in this module is primarily designed for generating datasets, training models, and performing analyses of the trained models. The input and output parameters, along with the intended usage of the four main functions, are outlined below:

\begin{function*}[\textbf{GenerateDataset}]Generate Dataset for Training:
$$\text{GenerateDataset }(str1, str2)$$ \par
\textbf{Input: }$str1$: filename for data; $str2$: filename for dataset;\par

\textbf{Description: }The file generated in Section \ref{sec:Program1} is in the Numpy array format and named ``$str1$". For efficient subsequent operations, this data is converted into a PyTorch tensor and structured as a dataset named ``$str2$". The program automatically shuffles the data, allocating $80\%$ as the training set and $20\%$ as the testing set.
\end{function*}

\begin{function*}[\textbf{LinearReg}] Linear Regression.
$$\text{LinearReg}(\mathbf{X},\mathbf{Y},\lambda)$$\par
\textbf{Input: }$\mathbf{X}\in\mathbb{R}^{N\times c}, \mathbf{Y}\in\mathbb{R}^{N}, \lambda\in\mathbb{R}$ \par
\textbf{Output: }$\beta\in\mathbb{R}^c, b\in\mathbb{R}$\par
\textbf{Description: }This function is used to solve a linear regression problem as outlined in Section \ref{sec:LinR}, with the hyperparameter of the regularization term set to $\lambda$. The $i$-th row of the matrix $\mathbf{X}$ represents $X_i$, and the $i$-th element of the vector $\mathbf{Y}$ represents $Y_i$. The same convention applies throughout the following discussion.

\end{function*}

\begin{function*}
    [\textbf{LinearCls}]SVM:
    $$\text{LinearCls}(\mathbf{X},\mathbf{Y},\lambda)$$\par
    \textbf{Input: }$\mathbf{X}\in\mathbb{R}^{N\times c}, \mathbf{Y}\in\mathbb{R}^{N}, \lambda\in\mathbb{R}$ \par
    \textbf{Output: }$\beta\in\mathbb{R}^c, b\in\mathbb{R}$\par
    \textbf{Description: }The function is used to solve SVM problem in Section \ref{sec:LinC}, with the hyperparameter of the regularization term set to $\lambda$.
\end{function*}

\begin{function*}[\textbf{NetReg}]Neural Network Regression:
$$\text{NetReg}(\mathbf{X},\mathbf{Y},n_L,n_H,\lambda)$$\par
\textbf{Input: }$\mathbf{X}\in\mathbb{R}^{N\times c}; \mathbf{Y}\in\mathbb{R}^{N}; n_L\in \mathbb{N}^+; n_H\in \mathbb{N}^+; \lambda\in\mathbb{R}$\par
\textbf{Output: }$\hat{f}:\text{torch.nn.module}$\par
\textbf{Description: }The symbol $\hat{f}$ denotes a trained fully connected network designed specifically for regression problems. The network comprises $N_L$ layers, each with $N_H$ neurons, as described in Section \ref{sec:NnR}. Weight decay is used as the regularization term, with the hyperparameter $\lambda$ controlling its strength.
\end{function*}

\begin{function*}
    [\textbf{NetCls}] Neural Network Classification:
    $$\text{NetCls}(\mathbf{X},\mathbf{Y},n_L,n_H,\lambda)$$\par
    \textbf{Input: }$\mathbf{X}\in\mathbb{R}^{N\times c}; \mathbf{Y}\in\mathbb{R}^{N}; n_L\in \mathbb{N}^+; n_H\in \mathbb{N}^+; \lambda\in\mathbb{R}$\par
    \textbf{Output: }$\hat{f}:\text{torch.nn.module}$\par
    \textbf{Description: }The symbol $\hat{f}$ denotes a trained fully connected network designed specifically for classification problems. The network comprises $N_L$ layers, each with $N_H$ neurons, as described in Section \ref{sec:NnC}. Weight decay is used as the regularization term, with the hyperparameter $\lambda$ controlling its strength.
\end{function*}

\begin{function*}[\textbf{NetGrad}] Sensitive Analysis:
$$\text{NetGrad}(\mathbf{X},\mathbf{Y},\hat{f})$$\par
\textbf{Input: }$\mathbf{X}\in\mathbb{R}^{N\times c}, \mathbf{Y}\in\mathbb{R}^{N},\hat{f}: \text{torch.NN.module}$\par
\textbf{Output: }$g\in\mathbb{R}^c:$\par
\textbf{Description: }This function is utilized to quantify the sensitivity of individual features in the trained network $\hat{f}$. The sensitivity is gauged by evaluating the mean of the absolute values of the derivatives of the loss function with respect to each feature, taken over the dataset $\{\mathbf{X},\mathbf{Y}\}$. Specifically, for the $j$-th feature, the sensitivity is calculated as:
$$g_{(j)}=\frac{1}{N} \sum_{i=1}^N |\frac{\partial \mathcal{L}(\hat{f}(X_i),Y_i)}{\partial X_{i,(j)}}|,$$
where $X_{i,(j)}$ denotes the $j$-th feature of the input example $X_i$, $N$ is the total number of examples, and $\mathcal{L}$ represents the loss function. The term $g_{(j)}$ delivers an aggregate measure of how sensitive the loss function is to variations in the $j$-th feature, thus quantifying the importance of that feature in the learned representation captured by $\hat{f}$.
\end{function*}

\section{Searching a dimension formula}\label{sec:SearchingDimFormula}

\subsection{Virtual Dimension} The journey towards the dimension formula of affine Deligne-Lusztig varieties has a long history. For the affine Grassmannian case, Rapoport proposed the dimension formula in \cite{Rapoport2002}, drawing inspiration from Chai's earlier work \cite{Chai2000} on the length function of chains of $\sigma$-conjugacy classes. This conjecture was ultimately validated by a series of researchers, primarily \cite{Viehmann2006, Goertz2006, Hamacher2015, Zhu16}.

In this paper, our attention is on the affine flag case, a significantly more challenging problem. G\"{o}rtz, Haines, Kottwitz, and Reuman proposed a conjectural formula for $\dim X_w(b)$ for most pairs $(w, b)$ in \cite{Goertz2010}. This was partly inspired by the aforementioned dimension formula for the affine Grassmannian case. This conjecture was verified by He in \cite{He2014_virtdim} and \cite{He2021_cordial}. However, our understanding of the remaining cases, which include many crucial applications to number theory and the Langlands program, remains limited. We aim to broaden this understanding using machine learning.

We revisit the concept of virtual dimension introduced by He in \cite{He2014_virtdim}. This was inspired by the conjecture of G\"{o}rtz, Haines, Kottwitz, and Reuman. Let $w \in W_a$ and express it as $w = x t^\mu y$ as in section \ref{sec:invariants}. Define $\eta(w) = yx$ and
$$d_w(b) = \frac{1}{2}(\ell(w) + \eta(w))-\<\nu_b,\rho\> - \frac{1}{2}\text{def}(b).$$
He demonstrated in \cite{He2014_virtdim} and \cite{He2021_cordial} the following result.
\begin{theorem}Suppose $X_w(b)\ne \emptyset$. Then
\begin{enumerate}
\item $\dim X_w(b)\le d_w(b)$.
\item If $w$ is $2$-regular and $\mu-\nu_b$ is ``sufficiently large'', then $\dim X_w(b)=d_w(b)$.
\end{enumerate}
\end{theorem}

The virtual dimension formula is a practical approximation of the real dimension. The discovery of this formula took experts considerable time, stretching from Rapoport's lectures in 1996 \cite{Rapoport2002} to the introduction of the virtual dimension formula in the most general setting by He in 2012 \cite{He2014_virtdim}. This section aims to illustrate how machine learning could help us rediscover this formula, accelerating the research process counterfactually.

The dimension of affine Deligne-Lusztig varieties depends on two parameters, $w$ and $b$. We note that the $b$-part of the virtual dimension formula is relatively simpler to uncover than the $w$-part, and we omit the details of the learning process for now. In this section, we concentrate on the group $\SL_5$, the case $b=1$, and randomly generated elements $w$. We investigate how machine learning can shed light on the correlation between $w$ and the dimension $\dim X_w(1)$. Our method does not rely on prior knowledge of the dimension formula in the Grassmannian case or the virtual dimension formula.

The selection of appropriate input features is crucial for machine learning. We work with the affine Weyl group of type $\tilde A_4$, where each element $w=t^\l u \in W_a$ is made up of the translation part $\l$ and the finite part $u$. We include both parts as input features: $\l$ as a vector, and $u$ as a permutation denoted by $u=[u_1, u_2, u_3, u_4, u_5]$.

For classical Deligne-Lusztig varieties $X_w$ \cite{Deligne1976}, it's known that $\dim X_w=\ell(w)$. Therefore, we anticipate that the dimension of affine Deligne-Lusztig varieties $X_w(1)$ is also related to $\ell(w)$. Consequently, we include the length function for both $w$ and $u$ as input features.

\subsection{Complexity test} 

To evaluate the complexity of the mapping, we utilize neural networks. Specifically, we consider an $n_L$-layer fully connected network with ReLU activation and $n_H$ hidden neurons. This network can be expressed as:
$$\hat{f}_\theta(X)=\beta^\top\rho(W{n_L}...W_2\rho(W_1X)...), \ \ \theta={W_1,W_2,...,W_{n_L},\beta},$$
where $\hat{f}_{\theta}(X)$ is the predicted output of the neural network for input $X\in\mathbb{R}^{c}$, and $\theta$ is the set of all trainable parameters, including weights:
$$W_1\in\mathbb{R}^{n_H\times c}, W_2\in\mathbb{R}^{n_H\times n_H}, ..., W{n_L}\in\mathbb{R}^{n_H\times n_H},\beta\in\mathbb{R}^{n_H}.$$
We use the Rectified Linear Unit (ReLU) activation function, defined as $\rho(x)=\max(0,x)$, for each hidden layer.

To quantify the accuracy of the optimized $\hat{f}_{\theta^*}$ in approximating $f$ after training, we use two metrics: accuracy and mean error. Since the dimension is always an integer, we round the inferred results $\hat{f}_{\theta^*}(X_i)$ to the nearest integer for accuracy. Specifically, we define accuracy as:
$$\text{Accuracy}=\frac{1}{N}\sum_i^N \delta(Y_i,round(\hat{f}_{\theta^*}(X_i))),$$
where $N$ is the number of samples, and $\delta(\cdot)$ is the indicator function that outputs $1$ if the arguments are equal and $0$ otherwise. The mean error is defined as:
$$\text{Mean Error} = \frac{1}{N}\sum_i^N |Y_i-\hat{f}_{\theta^*}(X_i)|.$$

\textbf{Dataset 1.}
We describe the first dataset used in our experiments. We randomly choose $5000$ elements $w$ from the set $W_a$ such that $\ell(w) < 30$ and $X_w(1) \ne \emptyset$. For each $w = t^{\l}u$, we express $\l = [\l_1,\l_2,\l_3,\l_4,\l_5]$ and $u = [u_1,u_2,u_3,u_4,u_5]$ as a permutation. We then compute the dimension $\dim X_w(1)$ for all these $w$.

\textbf{Experiment 1.}
In this experiment, we use Dataset 1 to train a neural network to predict the dimension of $X_w(1)$ for each $w=t^{\l}u$, where $\l$ and $u$ are defined as above. Specifically, the input vector for each $w$ is defined as $X = [\l_1,\l_2,\l_3,\l_4,\l_5, u_1,u_2,u_3,u_4,u_5,\ell(u),\ell(w)]$, and the corresponding output is $Y = \dim X_w(1)$. We obtain a dataset of $5000$ samples $(X_i,Y_i)_{i=1}^{5000}$ and train a neural network on it. The results of this experiment are summarized in Table \ref{Tabel:E_1}.

\begin{table}[htbp] 
\caption{Testing error of different neural networks for Dataset 1}
  \label{Table:nn_small}
  \centering
  \begin{tabular}{|c|c|c|c|}
    \hline
    \diagbox[width=5em]{$n_L$}{$n_H$}& 10 & 20 & 40 \\
        \hline
    1  & 0.53  &	0.53	& 0.52 \\
        \hline
    2 & 0.53	& 0.53 &	0.51 \\
    \hline
    3  & 0.52	&  0.51  &	0.51  \\
    \hline    
  \end{tabular}
\label{Tabel:E_1}
\end{table}

\textbf{Analysis.}
The results of Experiment 1 show that the error obtained is relatively small. Furthermore, we find that increasing the number of layers $n_L$ and hidden units $n_H$ in the network does not significantly enhance the accuracy of the neural network in predicting the dimension. These observations lead us to hypothesize that a linear model may be sufficient to approximate the mapping $f$ from the input data to the output dimension. This hypothesis will be investigated in the following subsection with the performance of a linear model.

\subsection{Linear model}
A linear function without a bias term can be expressed as $\hat{f}_\theta(X)=\beta^\top X =\sum_{i=1}^c \beta_{[i]}X_{[i]}$, where $\beta_{[i]}$ denotes the $i$-th element of the vector $\beta$ which represents the coefficient of the $i$-th term. For linear models, the mean error are used to evaluate the approximation of $\hat{f}_{\theta^*}$ to $f$.

\textbf{Experiment 2.}
We use Dataset 1 to train a linear model to predict the dimension of $X_w(1)$ for each $w=t^{\l}u$, where $\l$ and $u$ are as previously described. The mean error of the linear model is $0.65$.

\textbf{Analysis.}
Recall that the presentation $w = xt^\mu y$ is crucial for understanding the properties of the affine Weyl group element $w$, where $x,y\in W_0$ and $\mu \in X_\ast(T)$ is dominant. However, this presentation is not unique. Imposing a subtle restriction on $x$ or $y$ yields two uniquely determined presentations, but these may be different and thus are not canonical. Alternatively, one may restrict to those elements $w = xt^\mu y$ where $\mu$ is \emph{dominant regular}, i.e., $\mu_1>\cdots>\mu_5$. For such $w$, there are uniquely determined $x,y\in W_0$ with $w = xt^\mu y$.

We consider a subset of Dataset 1 consisting of those $w$ where $\mu$ is dominant regular, hypothesizing that the chosen features are more meaningful on this subset. This subset consists of $1037$ entries. A newly trained linear model on this subset attains a mean test error of $0.62$, which confirms our expectation. Hence, we consider analyzing the linear model under the dominant regular constraints by resampling $5000$ points to form Dataset 2.

\textbf{Dataset 2.}
Note that any element with a regular translation part can be written as $w =x t^{\mu}y$ where $x,y\in W_0$ and $\mu = [\mu_1,\ldots,\mu_5]$ is dominant regular, i.e., $\mu_1>\mu_2>\mu_3>\mu_4>\mu_5$. To investigate the neural network's performance on regular translation parts, we randomly choose $5000$ elements from the set of all such elements, where $\mu$ is dominant regular with $7>\mu_1>\mu_2>\mu_3>\mu_4>\mu_5>-7$, and $X_w(1)\ne\emptyset$. We compute $\dim X_w(1)$ for all these $w$.

\textbf{Experiment 3.}
We use Dataset 2. For each $w=xt^{\mu}y$, set \begin{align*}X &= [\d(x(\a_{12})), \d(x(\a_{13})),\ldots,\d(x(\a_{45})), \ell(x), \\&\qquad \mu_1,\mu_2,\mu_3,\mu_4,\mu_5, \d(y^{-1}(\a_{12})) ,\ldots,\d(y^{-1}(\a_{45})), \ell(y), \ell(w)]\end{align*} and $$Y=\dim X_w(1).$$

The input features are explained in detail as follows. We write $\a_{ij} = e_i - e_j$ for the root and
$$\delta(\a_{ij}) = \begin{cases}1,&i>j,\\0,&i<j,\end{cases}$$
for the indicator function of the negative roots. A linear combination of the permutation values $x(1),\ldots,x(n)$ would be very hard to interpret mathematically, which is why we use the $\delta$-values.

It is important to note for the neural network, the difference between presenting a permutation as $\delta(x(\a_{ij}))$ or $\delta(x^{-1}(\a_{ij}))$ is significant. While these two presentations are \emph{mathematically} equivalent, transitioning from one to the other is a non-linear procedure \cite[section 5.2]{Williamson2023}.
We use $x$ and $y^{-1}$, which have more direct mathematical interpretations than their inverses $x^{-1}$ and $y$. Specifically, $y^{-1}$ indicatesthe \emph{Weyl chamber} of $w$, whereas $x$ typically indicates the Weyl chamber of the inverse $w^{-1}$.

We obtain $(X_i,Y_i)_{i=1}^{5000}$. Upon applying the neural network to this dataset, the results of this experiment are summarized in Table \ref{Table:E_A}. The average error is found to be $0.65$.


\begin{table}[h]
\caption{Coefficient of Experiment 3}
  \centering
  \begin{tabular}{|c|c|}
    \hline
    Feature & Corresponding Coefficient(s) \\
 \hline
$\d(x(\a_{ij})), \d(y^{-1}(\a_{ij}))$ & $\left(\begin{smallmatrix}
\,&0.12& -0.04& -0.05& -0.24\\
&&0.14& -0.05& -0.04\\
&&&0.13& -0.02\\
&&&& 0.15\\
\,&&&&\end{smallmatrix}\right)$,  $\left(\begin{smallmatrix}
\,&0.10& -0.02& -0.08& -0.03\\
&&0.07& 0.00& -0.06\\
&&&0.02& 0.04\\
&&&&0.06\\
\,&&&&\end{smallmatrix}\right)$\\
            \hline
    $\ell(x)$  &0.10\\
            \hline
    $\mu_i$  & [0.13, -0.09, -0.02, 0.08, -0.10] \\
        \hline
    $\ell(y)$ &  0.10\\
    \hline
        $\ell(w)$ &  \textbf{0.52}\\
    \hline
    \end{tabular}
\label{Table:E_A}
\end{table}

\textbf{Analysis.} 
From the above table, we observe that the length of $w$ is the most significant feature, with a coefficient approximately equal to $1/2$. We thus deduct this potential leading term from the dimension for our subsequent experiments. Specifically, the output $Y$ will be $\dim X_w(1) - \frac{1}{2}\ell(w)$. It is noteworthy that $x$ and $y$ belong to the finite Weyl group, and their contribution to the dimension should be limited. Conversely, the range of $\mu$ is unbounded and it's anticipated that $\mu$ could contribute to the potential leading term of the linear approximation of the dimension. We therefore hypothesize that the contribution of $\mu$ in the linear model is already encapsulated in the term $\frac{1}{2}\ell(w)$ (i.e., the contribution of $\mu$ is given by $\<\mu, \rho\>$). Consequently, we will eliminate $\mu_i$ from $X$.

\textbf{Experiment 4.}
The analysis of the dimension of affine Deligne-Lusztig varieties, even for smaller rank groups such as $\SL_3$, suggests that different Weyl chambers may exhibit different patterns (cf. \cite[section~7]{Goertz2006} and \cite{Beazley2009}). In \cite{He2007}, He introduced the technique of partial conjugation, which, to some extent, reduces the problem to the dominant chamber. Therefore, we primarily focus on elements in the dominant Weyl chamber, i.e., elements of the form $w=t^{\mu}y$ where $\mu$ is dominant regular.

\textbf{Dataset 3.}
This dataset is defined similarly to Dataset 2. We randomly choose $5000$ elements $w =t^{\mu}y$ where $y\in W_0$ and $\mu = (\mu_1,\mu_2,\mu_3,\mu_4,\mu_5)$ is dominant regular with $9>\mu_1>\mu_2>\mu_3>\mu_4>\mu_5> - 9$, and $X_w(1)\ne\emptyset$.

For each $w = t^{\mu}y$, we define $$X = [ \d(y^{-1}(\a_{12})),\ldots, \d(y^{-1}(\a_{45})),\ell(y) ]$$ and $$Y=\dim X_w(1) - \frac{1}{2}\ell(w).$$
We derive $\{X_i,Y_i\}_{i=1}^{5000}$ and apply a linear model. We obtain $\hat{f}_{\theta^*}$ with an average error of $0.30$. The coefficients are listed in Table \ref{Table:E_B}.


\begin{table}[htbp]
\caption{Coefficient of Experiment 4}
  \centering
  \begin{tabular}{|c|c|}
    \hline
    Feature & Corresponding Coefficient(s) \\
        \hline
    $\d(y^{-1}(\a_{ij})) $  & $\left(\begin{smallmatrix}
\,&0.02& 0.05& 0.10& 0.13\\
&&-0.05& 0.07& 0.11\\
&&&-0.07& 0.05\\
&&&&0.00\\
\,&&&&\end{smallmatrix}\right)$\\   
        \hline   
    $\ell(y)$ &\textbf{0.42}\\
    \hline
    \end{tabular}
\label{Table:E_B}
\end{table}

Recall that $\ell(y) = \sum_{i<j}\delta(y^{-1}(\a_{i,j}))$. If we consider this linear dependence, we observe that the actual leading coefficient in the above linear model is $\ell(y)$ times
\begin{align*}
&0.42+(0.02+0.05+0.10+0.13-0.05+0.07+0.11-0.07+0.05+0.00)\cdot \frac 1{10}\\&=0.46. 
\end{align*}

\textbf{Analysis.} 
We could conjecture that $\frac 12\ell(y)$ is the leading term, with an \enquote{error term} that is one order of magnitude smaller. Before progressing with our experiments, we introduce some general strategies and terminology for the training and interpretation of linear models.

In the context of linear regression, the inclusion of a regularization term is essential when dealing with highly linearly correlated features. This situation can cause instability and unreliable estimates of the regression coefficients. Regularization involves adding a penalty term to the loss function, generally based on the magnitudes of the regression coefficients.

Practically, the $\ell_2$-norm is a commonly selected regularization term. The loss function can be represented as:
$$\min \sum_{i=1}^{N} (Y_i-\beta^\top X_i)^2+\lambda \lVert\beta\rVert_2^2.$$
The $\ell_2$-regularization in linear regression results in a unique optimal solution, providing model stability and avoiding multiple solutions. The $\ell_2$-regularization penalizes larger regression coefficients proportionally, leading to more stable models with reduced overfitting and improved generalization.

For instance, if we only use one feature, $\ell(y)$, to train the linear model, which does not require a regularization term, we can obtain a model with an average error of $0.30$ and a coefficient of $0.47$. This coefficient is much closer to $0.5$ because of the $\ell_2$-regularization's preference for small and average coefficients when dealing with linear dependence. Therefore, a common strategy is to first identify the most important features, discard less important features, and then retrain the model's parameters \cite{han2015deep, berglund2021machine}. This process allows for obtaining more accurate results and helps in mitigating the effects of highly correlated features.

Additionally, the $\ell_1$-norm is also a common metric used in regularization terms and fidelity terms. It can be expressed in the following form:  
$$\min \sum_{i=1}^{N} |Y_i-\beta^\top X_i|+\lambda \lVert\beta\rVert_1.$$
When used as a fidelity term, the $\ell_1$-norm penalizes absolute differences between the model predictions and the true targets. This makes it more robust to outliers compared to the $\ell_2$ norm, which is more sensitive to large errors. 

Also, as a regularization term, the $\ell_1$ norm induces sparsity in the model parameters, driving many parameters close to zero. This performs automatic feature selection, removing uninformative features and improving interpretability. In contrast, the $\ell_2$-norm does not induce sparsity, but instead diffuses weight across all parameters. Apply the $\ell_1$ model. We obtain $\hat{f}_{\theta^*}$ with average error of $0.18$. The coefficients are listed in Table \ref{Table:E_l1}.

\begin{table}[htbp]
\caption{Coefficient of $\ell_1$ model for Experiment 4}
  \centering
  \begin{tabular}{|c|c|}
    \hline
    Feature & Corresponding Coefficient(s) \\
        \hline
    $\d(y^{-1}(\a_{ij})) $  & $\left(\begin{smallmatrix}
\,&0.00&0.00&0.00&0.00\\
&&0.00&0.00&0.00\\
&&&0.00&0.00\\
&&&&0.00\\
\,&&&&\end{smallmatrix}\right)$\\   
        \hline   
    $\ell(y)$ &\textbf{0.50}\\
    \hline
    \end{tabular}
\label{Table:E_l1}
\end{table}



\textbf{Experiment 5.}
We aim to investigate how the findings from the previous sections extend to other Weyl chambers. Given $w = xt^{\mu}y$, our goal is to estimate
$$Y=\dim X_w(1)-\frac 12\ell(w) .$$
If $x=1$, we understand that the leading term should be $\frac 12 \ell(y)$. Generally speaking, there are a number of intuitive ways to combine $x$ and $y$, particularly in light of the previously mentioned partial conjugation method. These include the mutual products $xy$ and $yx$, as well as the Demazure products $y\ast x$ and $y \triangleleft x$ \cite{He2009}. Considering that the $\delta$-values only slightly contribute to the linear model in the dominant chamber, we exclude them in this experiment.

We utilize Dataset 2 (arbitrary Weyl chamber). For each $w = xt^{\mu}y$, we set $$X = [\ell(x), \ell(y),\ell(xy),\ell(yx),\ell(yx),\ell(y \triangleleft x)].$$
This gives us the pair $(X_i,Y_i)$ for $i=1$ to $5000$. Upon applying a linear model, we obtain $\hat{f}_{\theta^*}$ with a mean error of $0.13$. The coefficients are provided in Table \ref{Table:E_1_1}.

\begin{table}[htbp]
\caption{Coefficient of Experiment 5}
  \label{Table:E_1_1}
  \centering
  \begin{tabular}{|c|c|}
    \hline
    Feature & Corresponding Coefficient(s) \\
        \hline
    $\ell(x)$  & 0.02 \\
        \hline   
     $\ell(y)$  & -0.05\\
        \hline   
    $\ell(xy)$ & -0.02 \\
    \hline    
    $\ell(yx)$ & \textbf{0.46}\\
    \hline    
    $\ell(y*x)$ & 0.04\\
    \hline    
    $\ell(y\triangleleft x)$ & 0.04 \\    
    \hline    
  \end{tabular}
\end{table}

\textbf{Analysis.}
It is evident that $\ell(yx)$ is the most influential feature, yet interpreting the coefficient $0.46$ poses a mathematical challenge. Note that while the input features are \emph{linearly independent}, numerous mathematical relationships exist between them, as shown by the inequality \begin{align*}\left\lvert\ell(x)-\ell(y)\right\rvert\le \ell(y\triangleleft x) \le \ell(yx)\le \ell(y\ast x).\end{align*}
We can speculate that $\frac 12\ell(y x)$ could be a suitable candidate for the leading term, with the remaining terms being relatively small. This speculation leads us to the linear model
\begin{align*}
\dim X_w(1)\approx \frac 12\left(\ell(w) + \ell(yx)\right),
\end{align*}
which aligns with the previously mentioned virtual dimension $d_w(1)$.

\textbf{Experiment 6.}
We now scrutinize this linear model, which we have identified as representing the virtual dimension. For all pairs $(w,b)$ such that $w\in W_a$ satisfies $\ell(w)<30$ and $[b]\in B(\SL_5)$ satisfies $X_w(b)\neq\emptyset$, we examine the difference $d_w(b)-\dim X_w(b)$ between the virtual dimension and the actual dimension of $X_w(b)$. The number of such pairs is 3,119,946. The results are presented in Table~\ref{Table:F}.

\begin{table}[ht]
\caption{Number of pair $(w,b)$ with certain value of virt.\ dim minus dim.}\label{Table:F}
  \centering
  \begin{tabular}{|c|c|c|c|c|c| }
    \hline
     value of $d_w(b)-\dim X_w(b)$ & $0$ & $1$ & $2$ & $3$ & $4$ \\
        \hline
    Amount of pair & 2020909  & 922482  &166386  & 9885 & 284\\
    \hline
    \end{tabular}
\end{table}

\textbf{Analysis.}
We observe that the difference $d_w(b) - \dim X_w(b)$ is always non-negative, and equals zero in the majority of cases. Both of these observations are well-documented, and represent major accomplishments in the field \cite{He2014_virtdim, He2021_cordial}. Furthermore, we notice that this difference also seems to be upper-bounded by $4$, a surprisingly small value considering the number of pairs $(w,b)$ and the large dimensions involved. This leads us to conjecture that $d_w(b)-\dim X_w(b)$ always has a reasonably small upper bound independent of $(w,b)$. We will explore that question further in section~\ref{sec:lowerBound}.

\section{Searching for important features}\label{sec:importantFeatures}
In this section, we revisit the group $\SL_5$ without imposing any restrictions on $[b]\in B(\SL_5)$ or $w\in W_a$.

\subsection{Detailed Introduction to SVM method}

To ensure self-containment and enhance understanding of the experimental outcomes, we provide a detailed introduction to the Support Vector Machine (SVM) model. This widely-used machine learning algorithm is primarily employed for classification tasks and will be utilized in our experiments on the non-emptiness pattern and the condition of dimension equalling virtual dimension.

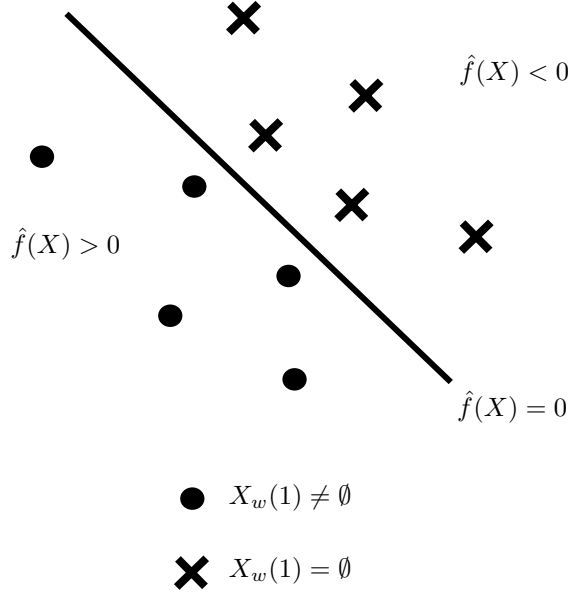
\begin{figure}[htb]
\centering
\begin{tikzpicture}[x=0.75pt,y=0.75pt,yscale=-1,xscale=1]

\draw  [draw opacity=0][fill={rgb, 255:red, 0; green, 0; blue, 0 }  ,fill opacity=1 ] (112,102.65) .. controls (112,99.53) and (114.71,97) .. (118.05,97) .. controls (121.39,97) and (124.1,99.53) .. (124.1,102.65) .. controls (124.1,105.77) and (121.39,108.3) .. (118.05,108.3) .. controls (114.71,108.3) and (112,105.77) .. (112,102.65) -- cycle ;
\draw [line width=2.25]    (54.25,15.75) -- (246,201) ;
\draw  [draw opacity=0][fill={rgb, 255:red, 0; green, 0; blue, 0 }  ,fill opacity=1 ] (159,147.65) .. controls (159,144.53) and (161.71,142) .. (165.05,142) .. controls (168.39,142) and (171.1,144.53) .. (171.1,147.65) .. controls (171.1,150.77) and (168.39,153.3) .. (165.05,153.3) .. controls (161.71,153.3) and (159,150.77) .. (159,147.65) -- cycle ;
\draw  [draw opacity=0][fill={rgb, 255:red, 0; green, 0; blue, 0 }  ,fill opacity=1 ] (36,87.65) .. controls (36,84.53) and (38.71,82) .. (42.05,82) .. controls (45.39,82) and (48.1,84.53) .. (48.1,87.65) .. controls (48.1,90.77) and (45.39,93.3) .. (42.05,93.3) .. controls (38.71,93.3) and (36,90.77) .. (36,87.65) -- cycle ;
\draw  [draw opacity=0][fill={rgb, 255:red, 0; green, 0; blue, 0 }  ,fill opacity=1 ] (100,167.65) .. controls (100,164.53) and (102.71,162) .. (106.05,162) .. controls (109.39,162) and (112.1,164.53) .. (112.1,167.65) .. controls (112.1,170.77) and (109.39,173.3) .. (106.05,173.3) .. controls (102.71,173.3) and (100,170.77) .. (100,167.65) -- cycle ;
\draw  [draw opacity=0][fill={rgb, 255:red, 0; green, 0; blue, 0 }  ,fill opacity=1 ] (162,199.65) .. controls (162,196.53) and (164.71,194) .. (168.05,194) .. controls (171.39,194) and (174.1,196.53) .. (174.1,199.65) .. controls (174.1,202.77) and (171.39,205.3) .. (168.05,205.3) .. controls (164.71,205.3) and (162,202.77) .. (162,199.65) -- cycle ;
\draw  [draw opacity=0][fill={rgb, 255:red, 0; green, 0; blue, 0 }  ,fill opacity=1 ] (111,259.65) .. controls (111,256.53) and (113.71,254) .. (117.05,254) .. controls (120.39,254) and (123.1,256.53) .. (123.1,259.65) .. controls (123.1,262.77) and (120.39,265.3) .. (117.05,265.3) .. controls (113.71,265.3) and (111,262.77) .. (111,259.65) -- cycle ;
\draw  [line width=3]  (146.87,69.99) -- (160.37,83.8)(160.99,69.69) -- (146.25,84.1) ;
\draw  [line width=3]  (189.87,104.99) -- (203.37,118.8)(203.99,104.69) -- (189.25,119.1) ;
\draw  [line width=3]  (135.87,10.99) -- (149.37,24.8)(149.99,10.69) -- (135.25,25.1) ;
\draw  [line width=3]  (196.87,49.99) -- (210.37,63.8)(210.99,49.69) -- (196.25,64.1) ;
\draw  [line width=3]  (251.87,120.99) -- (265.37,134.8)(265.99,120.69) -- (251.25,135.1) ;
\draw  [line width=3]  (109.87,289.99) -- (123.37,303.8)(123.99,289.69) -- (109.25,304.1) ;

\draw (249,35) node [anchor=north west][inner sep=0.75pt]   [align=left] {$\displaystyle \hat{f}( X) < 0$};
\draw (25,122) node [anchor=north west][inner sep=0.75pt]   [align=left] {$\displaystyle \hat{f}( X)  >0$};
\draw (248,204) node [anchor=north west][inner sep=0.75pt]   [align=left] {$\displaystyle \hat{f}( X) =0$};
\draw (134,251) node [anchor=north west][inner sep=0.75pt]   [align=left] {$\displaystyle X_{w}( 1) \neq \emptyset $};
\draw (134,288) node [anchor=north west][inner sep=0.75pt]   [align=left] {$\displaystyle X_{w}( 1) =\emptyset $};

\end{tikzpicture}

\caption{A demo for SVM}\label{fig:SVM}
\end{figure}

The primary objective of SVM is to identify an optimal hyperplane that effectively separates data points into different classes. In the case of binary classification, the hyperplane is chosen to maximize the margin, which refers to the distance between the hyperplane and the nearest data points from each class. In this context, we primarily focus on linear SVM for the sake of result interpretability. The equation of this hyperplane is given by:
$$\hat{f}(X)=\beta^\top  X - b = 0.$$
For instance, in the non-emptiness pattern experiments, SVM aims to establish a hyperplane that bifurcates the dataset into two regions. Specifically, on one side of the hyperplane, all instances satisfy $X_{w}(1) \neq \emptyset$; on the other side, all instances satisfy $X_{w}(1) = \emptyset$. This is depicted in Figure \ref{fig:SVM}.

When data is not inherently linearly separable, SVM often employs techniques to find an ``optimal" hyperplane. A commonly used evaluation criterion is the hinge loss, a margin-based loss function that penalizes misclassifications and encourages SVM to identify a decision boundary with a larger margin. For a binary classification problem, the hinge loss for a single data point is defined as:
$$\mathcal{L}_{\text{hinge}}(Y, \hat{f}(X)) = \max(0, 1 - Y \cdot \hat{f}(X)).$$
To optimize the SVM model, the aim is to minimize the sum of hinge losses across all training data points while incorporating a regularization term. This term helps prevent overfitting and controls the complexity of the learned model, often represented by the $\ell_2$-norm of the weight vector $\beta$. The SVM optimization problem can be formulated as:

$$\min  \sum_i \mathcal{L}_{\text{hinge}}(Y_i, \hat{f}(X_i))+\lambda \lvert \beta \rvert_2^2.$$
Solving this optimization problem allows SVM to learn a decision boundary that generalizes well to unseen data, leading to accurate classification or regression predictions.

During the inference stage, an input $X$ is classified as the first class if $\hat{f}(X)$ is greater than $0$, and classified as the second class otherwise. The coefficients represented by $\beta$ provide insightful information about the relationship between the input features and the class labels. Specifically, a positive $\beta_{(i)}$ suggests that as the value of $X_{(i)}$ increases, the likelihood of belonging to the first class also increases. Furthermore, a larger absolute value of $\beta_{(i)}$ signifies a stronger influence of the corresponding feature on the classification decision. Conversely, a negative $\beta_{(i)}$ indicates a negative relationship between the feature and the likelihood of belonging to the first class.

\subsection{Experiments on the non-emptiness pattern.\\}

\textbf{Data: }All $w=xt^{\mu}y\in W_a$ and $[b]\in B(\SL_5)$ respecting the conditions $\ell(w)<30$ and $\nu_b\le \mu$. The latter condition, known as Mazur's inequality, is a necessary condition for $X_w(b)\ne \emptyset$. (Dataset size: 8,705,879)

\textbf{Feature: }$$X = [  x_{ij},\mu_1,\dots,\mu_5, y^{-1}_{ij},\eta_{ij}=\d(\eta(w)(\a_{ij})), \ell(w), \nu_1,\dots,\nu_5,\l_1,\dots, \l_5]\in\mathbb{R}^{46}$$

\textbf{Output: }$Y = \begin{cases}
    1& \text{ if }X_w(1)\ne \emptyset\\
    -1& \text{ if }X_w(1)= \emptyset\\
\end{cases}$. 

\textbf{Result:}
Three models were utilized in our experiments, each executed $100$ times. The first model, SVM, yielded an average accuracy of $78.34\%$, with the average coefficients reported in Table \ref{Table:sen1_svm}. The second model, a single-layer neural ReLU classification network with $10$ neurons, achieved an average training accuracy of $87.14\%$ and an average test accuracy of $87.17\%$, with the average gradients documented in Table \ref{Table:sen1_1nn}. The third model, a three-layer neural ReLU classification network with $20$ neurons per layer, had an average training accuracy of $94.72\%$ and an average test accuracy of $94.74\%$, with the average gradients detailed in Table \ref{Table:sen1_3nn}.

\begin{table}[h]
\caption{Average Coefficient of SVM for non-empty }
\label{Table:sen1_svm}
  \centering
  \begin{tabular}{|c|c|}
    \hline
    Feature & Corresponding Coefficient(s) \\
        \hline
   $(x_{ij}),(y_{ij}^{-1}),(\eta_{ij})$   & $\left(\begin{smallmatrix}
\,&0.15&0.10&0.10&0.12\\
&&0.08&0.07&0.10\\
&&&0.08&0.10\\
&&&&0.15\\
\,&&&&\end{smallmatrix}\right)$,  $\left(\begin{smallmatrix}
\,&-0.08& -0.12& -0.15& -0.18\\
&&-0.07& -0.12& -0.15\\
&&&-0.07& -0.12\\
&&&&-0.08\\
\,&&&&\end{smallmatrix}\right)$,$\left(\begin{smallmatrix}
\,&-0.15& 0.10& \textbf{0.39}& \textbf{1.56}\\
&&-0.07& 0.00& \textbf{0.39}\\
&&&-0.07& 0.10\\
&&&&-0.15\\
\,&&&&\end{smallmatrix}\right)$ \\
\hline 
    $\mu_i$  & [0.01, 0.01, -0.00, -0.01, -0.02]\\
            \hline
    
    $\ell(w)$ & 0.07\\
    \hline
        $\nu_i$ & [\textbf{-1.14}, \textbf{-0.53}, 0.01, \textbf{0.55}, \textbf{1.16}]\\
    \hline
    $\l_i$ & [\textbf{0.68}, 0.30, 0.00, -0.30, \textbf{-0.67}]\\
    \hline
    \end{tabular}\\[2em]

\end{table}

\begin{table}[h]
\caption{Average Gradient of One-Layer NN for non-empty }
\label{Table:sen1_1nn}
  \centering
  \begin{tabular}{|c|c|}
    \hline
    Feature & Corresponding Gradient(s) \\
        \hline
 $(x_{ij}),(y_{ij}^{-1}),(\eta_{ij})$ & $\left(\begin{smallmatrix}
\,&0.07&0.05&0.04&0.05\\
&&0.04&0.03&0.04\\
&&&0.04&0.05\\
&&&&0.08\\
\,&&&&\end{smallmatrix}\right)$, $\left(\begin{smallmatrix}
\,&0.07&0.05&0.05&0.03\\
&&0.04&0.06&0.05\\
&&&0.04&0.06\\
&&&&0.08\\
\,&&&&\end{smallmatrix}\right)$, $\left(\begin{smallmatrix}
\,&0.08&\mathbf{0.21}&\mathbf{0.27}&\mathbf{0.69}\\
&&0.03&0.06&\mathbf{0.29}\\
&&&0.03&\mathbf{0.22}\\
&&&&0.08\\
\,&&&&\end{smallmatrix}\right)$\\
\hline
    $\mu_i$  & [0.08, 0.13, 0.14, 0.12, 0.07]\\
            \hline
    $\ell(w)$ & 0.05\\
    \hline
        $\nu_i$ & [\textbf{0.44}, \textbf{0.28}, 0.18, \textbf{0.28}, \textbf{0.45}]\\
    \hline
    $\l_i$ & [\textbf{0.23}, 0.14, 0.05, 0.13, \textbf{0.23}]\\
    \hline
    \end{tabular}
\end{table}


\begin{table}[h]
\caption{Average Gradient of Three-Layer NN for non-empty}
\label{Table:sen1_3nn}
  \centering
  \begin{tabular}{|c|c|}
    \hline
    Feature & Corresponding Gradient(s) \\
        \hline
$(x_{ij}),(y_{ij}^{-1}),(\eta_{ij})$ & $\left(\begin{smallmatrix}
\,&0.09&0.04&0.03&0.02\\
&&0.05&0.04&0.03\\
&&&0.05&0.04\\
&&&&0.09\\
\,&&&&\end{smallmatrix}\right)$, $\left(\begin{smallmatrix}
\,&0.09&0.04&0.03&0.03\\
&&0.06&0.04&0.03\\
&&&0.06&0.04\\
&&&&0.09\\
\,&&&&\end{smallmatrix}\right)$, $\left(\begin{smallmatrix}
\,&0.16&0.14&\mathbf{0.25}&\mathbf{0.30}\\
&&0.14&\mathbf{0.20}&\mathbf{0.26}\\
&&&0.14&0.14\\
&&&&0.16\\
\,&&&&\end{smallmatrix}\right)$\\
        \hline
    $\mu_i$  & [0.09, 0.18, \textbf{0.21}, 0.19, 0.09]\\
            \hline
    $\ell(w)$ & 0.06\\
    \hline
        $\nu_i$ & [0.14, 0.10, 0.09, 0.11, 0.14]\\
    \hline
    $\l_i$ & [0.17, 0.16, 0.16, 0.16, 0.18]\\
    \hline
    \end{tabular}
\end{table}

\textbf{Analysis.}
It is evident that the \enquote{hook}-part of $\eta$ encapsulates important features, specifically $\eta_{12}, \eta_{13}, \eta_{14}$, $\eta_{15},\eta_{25},\eta_{35}, \eta_{45}$. The No Levi Obstruction (NLO) criterion for non-emptyness of ADLV, developed in \cite{Goertz2010} and further refined in \cite{Goertz2015} \cite{Viehmann2021}, suggests that the support of $\eta(w)$ plays a crucial role in the non-empty pattern. This connects to the values of $\eta_{ij}$, though the exact relationship is intricate and non-linear, representing an interesting problem in itself.

As seen from Table \ref{Table:sen1_svm}, another important feature is the Newton point $\nu_b$ together with its best integral approximation $\lambda(b)$. Several known results indicate that $X_w(b)=\emptyset$ holds frequently if $[b]$ is very large and only occasionally if $[b]$ is very small \cite{He2021_cordial, Schremmer2022}.

The three layer neural network achieves a high accuracy, but the gradient is hard to interpret mathematically. This is partly due to inherent properties of the machine learning model (cf.\ section~\ref{sec:caveatsChoiceOfModel}), as well as the underlying problem: Since the target function $f$ is discrete, the gradient of the smooth approximation $\hat f$ at individual points is hard to interpret. Nonetheless, this lack of explanation can inspire future mathematical research, directing our attention to get a more precise understanding of the underlying mathematical problem.

\

\subsection{Experiments on the dimension.\\}

\textbf{Data: }All $w=xt^{\mu}y\in W_a$ and $[b]\in B(\SL_5)$, with the conditions $\ell(w)<30$ and $X_w(b)\ne \emptyset$. (Dataset size: 3,119,946)

\textbf{Feature: }$$X = [x_{ij},\mu_1,\dots,\mu_5, y^{-1}_{ij},\eta_{ij} , \ell(w), \nu_1,\dots,\nu_5,\l_1,\dots, \l_5]\in\mathbb{R}^{46}$$

\textbf{Output: }$Y = \dim X_w(b)$.

\textbf{Result:}
Two models were utilized in our experiments, each executed $100$ times. The first model, a single-layer neural ReLU classification network with $10$ neurons, achieved an average training accuracy of $77.02\%$ and an average training error of $0.33$, an average test accuracy of $77.06\%$ and an average test error of $0.33$. The average gradients can be found in Table \ref{Table:sen2_1nn}. The second model, a three-layer neural ReLU classification network with $20$ neurons per layer, reached an average training accuracy of $92.12\%$ and an average training error of $0.16$, an average test accuracy of $92.10\%$ and an average test error of $0.16$. The average gradients are detailed in Table \ref{Table:sen2_3nn}.

\begin{table}[h]
\caption{Average Gradient of One-Layer NN for dimension }
\label{Table:sen2_1nn}
  \centering
  \begin{tabular}{|c|c|}
    \hline
    Feature & Corresponding Gradient(s) \\
        \hline
  $(x_{ij}),(y_{ij}^{-1}),(\eta_{ij})$   & $\left(\begin{smallmatrix}
\,&0.22&0.13&0.05&0.04\\
&&0.26&0.13&0.05\\
&&&0.26&-0.13\\
&&&&0.22\\
\,&&&&\end{smallmatrix}\right)$,  $\left(\begin{smallmatrix}
\,&0.18&0.16&0.0.08&0.03\\
&&0.20&0.18&0.08\\
&&&0.20&0.16\\
&&&&0.18\\
\,&&&&\end{smallmatrix}\right)$,$\left(\begin{smallmatrix}
\,&0.21&0.26&0.32&0.37\\
&&0.18&0.25&0.31\\
&&&0.18&0.26\\
&&&&0.21\\
\,&&&&\end{smallmatrix}\right)$ \\
\hline 
        
    $\mu_i$  & [0.34, 0.20, 0.06, 0.20, 0.35]\\
            \hline
    $\ell(w)$ & 0.46\\
    \hline
        $\nu_i$ & [\textbf{4.19}, \textbf{2.08}, 0.05, \textbf{2.07}, \textbf{4.19}]\\
    \hline
    $\l_i$ & [\textbf{2.01}, \textbf{1.00}, 0.05, \textbf{1.01}, \textbf{2.02}]\\
    \hline
    \end{tabular}
\end{table}

\begin{table}[h]
\caption{Average Gradient of Three-Layer NN for dimension}
\label{Table:sen2_3nn}
  \centering
  \begin{tabular}{|c|c|}
    \hline
    Feature & Corresponding Gradient(s) \\
        \hline
   $(x_{ij}),(y_{ij}^{-1}),(\eta_{ij})$   & $\left(\begin{smallmatrix}
\,&0.18&0.12&0.08&0.08\\
&&0.24&0.13&0.08\\
&&&0.24&0.12\\
&&&&0.18\\
\,&&&&\end{smallmatrix}\right)$,  $\left(\begin{smallmatrix}
\,&0.16&0.12&0.08&0.18\\
&&0.21&0.13&0.08\\
&&&0.22&0.11\\
&&&&0.16\\
\,&&&&\end{smallmatrix}\right)$,$\left(\begin{smallmatrix}
\,&0.40&0.41&0.45&0.47\\
&&0.42&0.42&0.46\\
&&&0.41&0.41\\
&&&&0.40\\
\,&&&&\end{smallmatrix}\right)$ \\
\hline       
    $\mu_i$  & [0.37, 0.32, 0.31, 0.31, 0.37]\\
            \hline
    $\ell(w)$ & 0.47\\
    \hline
        $\nu_i$ & [\textbf{4.22}, \textbf{2.03}, 0.29, \textbf{2.03}, \textbf{4.22}]\\
    \hline
    $\l_i$ & [\textbf{2.01}, \textbf{1.02}, 0.14, \textbf{1.00}, \textbf{2.00}]\\
    \hline
    \end{tabular}
\end{table}

\textbf{Analysis.}
It can be observed that the average gradients for these neural networks closely align with the gradient of the linear model discussed in \S 4. In fact, the dimension equals the virtual dimension for $64.8\%$ of the dataset, and for most of the remaining data points, the difference is just $1$. The gradient of this linear function seems to dominate the more nuanced behavior of the neural network, resulting in the $92.1\%$ accuracy. To obtain a more insightful average gradient, a comparison between the dimension and the virtual dimension should be considered.

\

\subsection{Experiments on the condition virtual dim.\ = dim.\\}

\textbf{Data: }All $w=xt^{\mu}y\in W_a$ and $[b]\in B(\SL_5)$, with the conditions $\ell(w)<30$ and $X_w(b)\ne \emptyset$. (Dataset size: 3,119,946)

\textbf{Feature: }$$X = [ x_{ij},\mu_1,\dots,\mu_5, y^{-1}_{ij},\eta_{ij}=\d(\eta(w)(\a_{ij})), \ell(w), \nu_1,\dots,\nu_5,\l_1,\dots, \l_5]\in\mathbb{R}^{46}$$

\textbf{Output: }$Y = \begin{cases}
    1& \text{ if } VD \ne \dim X_w(b)\\
    -1& \text{ if } VD = \dim X_w(b)\\
\end{cases}$. 

\textbf{Result:}
Three models were utilized in our experiments, each executed $100$ times. The first model, an SVM, achieved an average accuracy of $83.13\%$, with the average coefficients shown in Table \ref{Table:sen_VDc_svm}. The second model, a single-layer neural ReLU classification network with 10 neurons, recorded an average training accuracy of $87.52\%$ and an average test accuracy of $87.55\%$, with the average gradients shown in Table \ref{Table:sen_VDc_1nn}. The third model, a three-layer neural ReLU classification network with $20$ neurons per layer, attained an average training accuracy of $96.66\%$ and an average test accuracy of $96.67\%$, with the average gradients presented in Table \ref{Table:sen_VDc_3nn}.

\begin{table}[h]
\caption{Average Coefficient of SVM for VD=Dim}
\label{Table:sen_VDc_svm}
  \centering
  \begin{tabular}{|c|c|}
    \hline
    Feature & Corresponding Coefficient(s) \\
        \hline
   $(x_{ij}),(y_{ij}^{-1}),(\eta_{ij})$   & $\left(\begin{smallmatrix}
\,&-0.55&-0.10&0.21&0.36\\
&&-0.87& -0.13& 0.22\\
&&&-0.87 & -0.10\\
&&&&-0.55\\
\,&&&&\end{smallmatrix}\right)$,  $\left(\begin{smallmatrix}
\,&0.39& 0.28 & -0.02 & -0.27\\
&&0.67&0.32& -0.02\\
&&& 0.67& 0.27\\
&&&&0.40\\
\,&&&&\end{smallmatrix}\right)$,$\left(\begin{smallmatrix}
\,&\textbf{1.18} & 0.91 & 0.70 & 0.73\\
&&\textbf{1.57}& 0.95& 0.70,\\
&&&\textbf{1.57}& 0.91\\
&&&&\textbf{1.18}\\
\,&&&&\end{smallmatrix}\right)$ \\
\hline       
    $\mu_i$  & [0.02, -0.14, -0.00, 0.15, -0.02]\\
            \hline
    $\ell(w)$ & -0.19\\
    \hline
        $\nu_i$ & [\textbf{1.02}, 0.44, 0.00, -0.45, \textbf{-1.02}]\\
    \hline
    $\l_i$ & [-0.24, -0.10, 0.00, 0.10, 0.24]\\
    \hline
    \end{tabular}
\end{table}

\begin{table}[h]
\caption{Average Gradient of One-Layer NN for VD=Dim}
\label{Table:sen_VDc_1nn}
  \centering
  \begin{tabular}{|c|c|}
    \hline
    Feature & Corresponding Gradient(s) \\
        \hline
           $(x_{ij}),(y_{ij}^{-1}),(\eta_{ij})$   & $\left(\begin{smallmatrix}
\,&0.17& 0.06& 0.07& 0.11\\
&&\textbf{0.38}& 0.07& 0.07\\
&&&\textbf{0.44}& 0.06\\
&&&&0.18\\
\,&&&&\end{smallmatrix}\right)$,  $\left(\begin{smallmatrix}
\,&0.12& 0.07& 0.04& 0.09\\
&&0.21& 0.08& 0.04\\
&&&0.21& 0.08\\
&&&&0.13\\
\,&&&&\end{smallmatrix}\right)$,$\left(\begin{smallmatrix}
\,&\textbf{0.38}& 0.32& 0.22& 0.22\\
&&\textbf{0.57}& 0.30& 0.23\\
&&&\textbf{0.56}& 0.32\\
&&&&\textbf{0.39}\\
\,&&&&\end{smallmatrix}\right)$ \\
\hline 
    $\mu_i$  & [0.06, 0.16, 0.24, 0.16, 0.06]\\
            \hline
    $\ell(w)$ & 0.06\\
    \hline
        $\nu_i$ & [\textbf{0.37}, 0.18, 0.18, 0.18, \textbf{0.38}]\\
    \hline
    $\l_i$ & [0.12, 0.06, 0.06, 0.06, 0.12]\\
    \hline
    \end{tabular}
\end{table}

\begin{table}[h]
\caption{Average Gradient of Three-Layer NN for VD=Dim}
\label{Table:sen_VDc_3nn}
  \centering
  \begin{tabular}{|c|c|}
    \hline
    Feature & Corresponding Gradient(s) \\
        \hline
   $(x_{ij}),(y_{ij}^{-1}),(\eta_{ij})$   & $\left(\begin{smallmatrix}
\,&0.13 & 0.06 & 0.04 & 0.06\\
&&\textbf{0.35}& 0.07& 0.04\\
&&&\textbf{0.31}& 0.06\\
&&&&0.12\\
\,&&&&\end{smallmatrix}\right)$,  $\left(\begin{smallmatrix}
\,&0.10& 0.04& 0.04& 0.05\\
&&0.20& 0.05& 0.04\\
&&&0.19& 0.04\\
&&&&0.10\\
\,&&&&\end{smallmatrix}\right)$,$\left(\begin{smallmatrix}
\,&0.27& 0.27& 0.24& 0.22\\
&&\textbf{0.39}& 0.25& 0.24\\
&&&\textbf{0.38}& 0.28\\
&&&&0.25\\
\,&&&&\end{smallmatrix}\right)$ \\
\hline     
    $\mu_i$  & [0.08, 0.23, \textbf{0.34}, 0.22, 0.07]\\
            \hline
    $\ell(w)$ & 0.05\\
    \hline
        $\nu_i$ & [\textbf{0.31}, 0.21, 0.20, 0.20, \textbf{0.31}]\\
    \hline
    $\l_i$ & [0.12, 0.07, 0.07, 0.07, 0.12]\\
    \hline
    \end{tabular}
\end{table}

\textbf{Analysis.}
We see that $\eta_{i,j}$ are important features. We know that if $\eta(w)$ is small, for example, if $\eta(w)$ is a partial Coxeter element, then the dimension equals virtual dimension \cite{He2022}. If $\eta(w)$ is large and close to the longest element, then dimension is unlikely to be equal to virtual dimension.

Moreover, it is known if $y=1$ or, under additional hypotheses, if $x=w_0$, then dimension must also be equal to virtual dimension \cite{Milicevic2020}. This explains the signs of the $x_{ij}$ and $y_{ij}$ in Table \ref{Table:sen_VDc_svm}.

It is known in general, cf.\ \eqref{eq:purityConsequence}, that the difference $d_w(b)-\dim X_w(b)$ is maximal for large $[b]$. This explains why the Newton point is an important feature.


Table~\ref{tab:probabilityDimVirtDimFailure} gives the proportion of elements $w$, grouped by the length of $\eta(w)$, where the virtual dimension is not equal to dimension for some $b$ (i.e.\ the non-cordial elements in the sense of \cite{Milicevic2020}).

\begin{table}[ht]\small
\caption{number of elements with dim=virdim fail}\label{tab:probabilityDimVirtDimFailure}
\label{Table:Delta}
  \centering
  \begin{tabular}{|c|c|c|c|c|c|c|c|c|c|c|c|}
    \hline
    $\ell(\eta(w))$  & $0$ & $1$ & $2$ & $3$ & $4$& $5$ & $6$&$7$&$8$&$9$&$10$ \\
        \hline
    $\#$ non-cordial  & 0  & 0  &0  & 2696 & 8316  & 18232 &18152 &17651&10039&5284&1175\\    
        \hline
    $\#$ cordial  &  1271 & 5742  & 11191 & 21255 & 24754  & 31172 &24780 &21292&11155&5664&1225\\    
        \hline    
    Proportion & 0&0&0&11\% & 25\% & 37\%& 42\% & 45\% & 47\%&48\%& 49\% \\
    \hline
    \end{tabular}
\end{table}

\subsection{Statistics of the difference of virtual dim.\ - dim.\\} The experiments and analysis above indicate that the virtual dimension is a good approximation of the dimension for nonempty $X_w(b)$. A natural question is to further study the difference between the dimension and virtual dimension for nonempty $X_w(b)$.

In this part, we will provide a numerical analysis of the the difference between the dimension and virtual dimension $\Delta_w(b)=d_w(b)-\dim X_w(b)$ for nonempty $X_w(b)$ of dataset 2, and hope that the analysis will guide us to get the pattern of difference.

\begin{table}[ht]
\caption{$\Delta_w(b)$}
\label{Table:Delta2}
  \centering
  \begin{tabular}{|c|c|c|c|c|c|c|}
    \hline
      & $\Delta_w(b)=0$ & $\Delta_w(b)=1$ & $\Delta_w(b)=2$ & $\Delta_w(b)=3$ & $\Delta_w(b)=4$& All\\
        \hline
    Amount of data  & 2020909 & 922482 &166386 &9885 & 284 &3119946 \\    

        \hline
    Proportion of data &65\% & 30\%& 5\%&0.3\% &  0.01\% &100\%\\
    \hline
    \end{tabular}
\end{table}

As exhibited in Table \ref{Table:Delta2}, within the dataset, the virtual dimension represents the theoretical upper bound of the actual dimension. The maximum achievable accuracy rate was $65.82\%$. Furthermore, the predominant $\Delta_w(b)$ were of magnitude $1$.  

It seems that the percentage of pairs $(w,b)$ with $\Delta_w(b)$ larger than a given bound decreases rapidly. We further expect that $\Delta_w(b)$ might be bounded from above by a constant depending only on $n$, i.e.\ the group $\SL_n$. We will investigate this question by mathematical methods in section~\ref{sec:lowerBound}.

\subsection{Experiments on the irreducible components.\\}

In this section, we investigate the number of top dimensional irreducible components of $X_w(b)$ up to the $\bJ_b$-action. The analogous question has been solved for affine Deligne--Lusztig varieties in the affine Grassmannian \cite{Zhou2020, Nie2018}. Here, the key ingredient is the dimension of the \emph{weight space of the highest-weight Verma module} $V_\mu(\lambda)$. This is a commonly studied object in the representation theory of Lie algebras, and we refer the reader to any of the corresponding textbooks for the exact definition. For now, we remark that this dimension is a positive integer that is computable in terms of $\mu \in X_\ast(T)$ (which depends only on $w$) and $\lambda(b)\in X_\ast(T)$ (which only depends on $b$). Our first experiment uses the data set from the previous experiment, restricted to those ADLV whose dimension agrees with virtual dimension.

(1) \textbf{Data: }All $w=xt^{\mu}y\in W_a$ and $[b]\in B(\SL_5)$ with $\ell(w)<30$ and $\dim X_w(b)= d_w(b)$. (Dataset size: 2,020,909)

\textbf{Feature: }$$X = [  x_{ij},\mu_1,\dots,\mu_5, y^{-1}_{ij},\eta_{ij} , \ell(w), \nu_1,\dots,\nu_5,\l_1,\dots, \l_5, \text{dim} V_{\mu}(\lambda)  ]\in\mathbb{R}^{47}$$

\textbf{Output: }$Y = \sharp \bJ_b \backslash \Sigma^{\text{top}} X_w(b)$. 

\textbf{Result:}
Two models were utilized in our experiments, each executed $100$ times. 
The first model, a single-layer neural ReLU classification network with $10$ neurons, achieved an average
training accuracy of $67.15\%$ and an average training error of $0.48$, an average test accuracy of $67.06\%$ and an average test error of $0.48$.  The average gradients can be found in Table \ref{Table:6621}. The second model, a three-layer neural ReLU classification network with $20$ neurons per layer, reached an average training accuracy of $75.96\%$ and an average training error of $0.33$, an average test accuracy of $75.92\%$ and an average test error of $0.33$. The average gradients are shown in Table \ref{Table:6622}.


\begin{table}[h]
\caption{Average Gradient of One-Layer NN for irreducible components }
\label{Table:6621}
  \centering
  \begin{tabular}{|c|c|}
    \hline
    Feature & Corresponding Gradient(s) \\
\hline
   $(x_{ij}),(y_{ij}^{-1}),(\eta_{ij})$   & $\left(\begin{smallmatrix}
\,&0.23& 0.13& 0.07& 0.21\\
&&0.29& 0.12& 0.07\\
&&&0.29& 0.13\\
&&&&0.23\\
\,&&&&\end{smallmatrix}\right)$,  $\left(\begin{smallmatrix}
\,&0.25& 0.14& 0.04& 0.19\\
&&0.29& 0.15& 0.04\\
&&&0.29& 0.14\\
&&&&0.25\\
\,&&&&\end{smallmatrix}\right)$,$\left(\begin{smallmatrix}
\,&0.34& 0.26& 0.35& \textbf{0.62}\\
&&0.41& \textbf{0.54}& 0.35\\
&&& 0.43& 0.26\\
&&&&0.34\\
\,&&&&\end{smallmatrix}\right)$ \\
\hline        
    $\mu_i$  & [0.04, 0.11, 0.06, 0.11, 0.04]\\
            \hline
    $\ell(w)$ & 0.08\\
    \hline
        $\nu_i$ & [\textbf{0.56}, 0.30, 0.03, 0.30, \textbf{0.55}]\\
    \hline
    $\l_i$ & [0.25, 0.13, 0.02, 0.13, 0.25]\\
    \hline
    $\text{dim}V_{\mu}(\l)$ & 0.02\\
    \hline
    \end{tabular}
\end{table}

\begin{table}[h]
\caption{Average Gradient of Three-Layer NN for irreducible components}
\label{Table:6622}
  \centering
  \begin{tabular}{|c|c|}
    \hline
    Feature & Corresponding Gradient(s) \\
        \hline
   $(x_{ij}),(y_{ij}^{-1}),(\eta_{ij})$   & $\left(\begin{smallmatrix}
\,&0.25& 0.16& 0.10& 0.17\\
&&0.31& 0.15& 0.10\\
&&&0.31& 0.16\\
&&&&0.25\\
\,&&&&\end{smallmatrix}\right)$,  $\left(\begin{smallmatrix}
\,&0.25& 0.14& 0.08& 0.16\\
&&0.30& 0.16& 0.08\\
&&&0.31& 0.14\\
&&&&0.26\\
\,&&&&\end{smallmatrix}\right)$,$\left(\begin{smallmatrix}
\,&0.54& 0.43& 0.63& 0.55\\
&&\textbf{0.69}& \textbf{0.87}& 0.62\\
&&&\textbf{0.72}& 0.42\\
&&&&0.51\\
\,&&&&\end{smallmatrix}\right)$ \\
\hline       
    $\mu_i$  & [0.08, 0.28, 0.33, 0.27, 0.08]\\
            \hline
    $\ell(w)$ & 0.07\\
    \hline
        $\nu_i$ & [0.50, 0.30, 0.15, 0.30, 0.50]\\
    \hline
    $\l_i$ & [0.24, 0.13, 0.06, 0.13, 0.24]\\
    \hline
    $\text{dim}V_{\mu}(\l)$&0.03\\
    \hline
    \end{tabular}
\end{table}

\textbf{Analysis.}
We observe that the variables $\nu_1,\nu_5,\eta_{15},\eta_{23},\eta_{24},\eta_{34},\mu_3$ play a sensitive role in the approximated function $\hat{f}$. Interestingly, the simpler model tends to overlook $\eta_{23},\eta_{24},\eta_{34},\mu_3$, while these variables appear to be crucial for the more complex model. This suggests that the function $f$ exhibits a complex relationship with respect to $\eta_{23},\eta_{24},\eta_{34},\mu_3$. Moreover, the accuracy achieved in this experiment is lower than that of the previous one, which focused on the dimension problem. This suggests that the problem of determining irreducible components presents greater complexity.

We remark that the feature $\dim V_\mu(\lambda)$ does not seem to be particularly influential. Repeating the experiment without this feature, we retain almost the same accuracy. In order to obtain further insight into the most well-behaved situations, we restrict our attention to a certain subset of elements, which are known to enjoy the most convenient properties.

 (2) \textbf{Data: }All $w=xt^{\mu}y\in W_a$ where $y=1$ and $[b]\in B(\SL_5)$ with $\ell(w)<30$ and $X_w(b)\neq\emptyset$. In this case, it is known that $\dim X_w(b)= d_w(b)$. (Dataset size: 43,986)

\textbf{Feature: }$$X = [  x_{ij},\mu_1,\dots,\mu_5, \ell(w), \nu_1,\dots,\nu_5,\l_1,\dots, \l_5,\text{dim} V_{\mu}(\lambda) ]\in\mathbb{R}^{47}$$

\textbf{Output: }$Y = \sharp \bJ_b \backslash \Sigma^{\text{top}} X_w(b)$.

\textbf{Result: }Two models were utilized in our experiments, each executed $100$ times. The first model, a single-layer neural ReLU classification network with $10$ neurons, achieved an average training accuracy of $67.82\%$ and an average training error of $0.49$, an average test accuracy of $67.67\%$ and an average test error of $0.49$. The average gradients can be found in Table \ref{Table:6631}. The second model, a three-layer neural ReLU classification network with $20$ neurons per layer, reached an average training accuracy of $81.80\%$ and an average training error of $0.26$, an average test accuracy of $81.60\%$ and an average test error of $0.27$. The average gradients are shown in Table \ref{Table:6632}.

\begin{table}[h]
\caption{Average Gradient of One-Layer NN for irreducible components }
\label{Table:6631}
  \centering
  \begin{tabular}{|c|c|}
    \hline
    Feature & Corresponding Gradient(s) \\
        \hline
    $(x_{ij})$  & $\left(\begin{smallmatrix}
\,&0.21& 0.20& 0.25& \textbf{0.29}\\
&&0.16& \textbf{0.27}& 0.25\\
&&&0.17& 0.20\\
&&&&0.20\\
\,&&&&\end{smallmatrix}\right)$\\
            \hline
    $\mu_i$  & [0.11, 0.13, 0.09, 0.12, 0.10]\\
            \hline
  
    $\ell(w)$ & 0.04 \\
    \hline
        $\nu_i$ & [0.12, 0.10, 0.05, 0.10, 0.11]\\
    \hline
    $\l_i$ & [0.07, 0.05, 0.03, 0.04, 0.07]\\
    \hline
        $\text{dim}V_{\mu}(\l)$&  \textbf{0.29}\\
    \hline
    \end{tabular}
\end{table}

\begin{table}[h]
\caption{Average Gradient of Three-Layer NN for irreducible components}
\label{Table:6632}
  \centering
  \begin{tabular}{|c|c|}
    \hline
    Feature & Corresponding Gradient(s) \\
        \hline
        $(x_{ij})$  & $\left(\begin{smallmatrix}
\,&0.21& 0.20& 0.23& 0.20\\
&&0.21& \textbf{0.26}& 0.23\\
&&&0.21& 0.20\\
&&&&0.21\\
\,&&&&\end{smallmatrix}\right)$\\
            \hline
    $\mu_i$  & [0.10, 0.19, \textbf{0.25}, 0.20, 0.11]\\
            \hline

    $\ell(w)$ & 0.02 \\
    \hline
        $\nu_i$ & [0.09, 0.09, 0.13, 0.09, 0.09]\\
    \hline
    $\l_i$ & [0.05, 0.04, 0.05, 0.04, 0.05]\\
    \hline
            $\text{dim}V_{\mu}(\l)$&  \textbf{0.28}\\
    \hline
    \end{tabular}
\end{table}

\textbf{Analysis.}
The restriction to $y=1$ implies that $X_w(b)$ is equidimensional $\dim X_w(b) = d_w(b)$ whenever $X_w(b)\neq\emptyset$ \cite{Milicevic2021}. In other words, all irreducible components are top dimensional.

We see that the accuracy of the single-layer neural network does not change much compared to the previous experiment. However, the gradients are vastly different. The biggest contribution in Table~\ref{Table:6631} comes from the dimension of the weight space $\dim V_\mu(\lambda)$, which only made a tiny contribution in the previous experiment. In the case of three-layer neural networks, the restriction to $y=1$ brings a substantial improvement in accuracy, and again the contribution of $\dim V_\mu(\lambda)$ becomes significantly larger.

If $x =w_0$ is the longest element of the Weyl group, we know that the irreducible components of $X_w(b)$ correspond one-to-one to the irreducible components of the affine Deligne-Lusztig variety $X_\mu(b)$ inside the affine Grassmannian (following the proof of \cite[Theorem~10.1]{He2014_virtdim}). For the latter kind of affine Deligne-Lusztig varieties, the number of $\bJ_b$-orbits of irreducible components has been predicted by Chen-Zhu, and their conjecture has been fully established by Zhou-Zhu \cite{Zhou2020} as well as Nie \cite{Nie2018}. In this case, we know that $Y = \dim V_\mu(\lambda)$.

More generally, the same conclusion holds whenever $x$ is the longest element of a Levi subgroup, following the Hodge-Newton decomposition method of Görtz-He-Nie \cite{Goertz2015}. For general $x$, we may expect that the number of irreducible components is much smaller. Nonetheless, it should not be a surprise that $\dim V_\mu(\lambda)$ is the input feature with the highest overall contribution, as measured by the average gradient. It is not quite clear why this was not the case in the previous experiment, since these two cases are related by the partial conjugation method, which is independent of $\mu$. This could be an artifact of our limited data set. However, the situation remains overall mysterious, and we invite the interested readers to further explore this phenomenon through mathematical insight or ML assisted research.

Overall, we may summarize that the problem of enumerating irreducible components allows for a fairly accurate solutions using single-layer or three-layer neural networks. This gives hope that further mathematical progress on this problem should be possible. Moreover, the second subset does indeed seem to be better behaved for studying this problem.

\section{Lower bound on the dimension} \label{sec:lowerBound}

In section~\ref{sec:importantFeatures}, we developed machine
learning models that not only enable us to recover previously known
results, but also lead to new conjectures and research questions. In this section, we study the question on the lower bound of the dimension of ADLV whenever it is non-empty. In other words, we study the upper bound of the difference $d_w(b) - \dim X_w(b)$ whenever $X_w(b)\neq\emptyset$. For $\bG = \SL_n$, we will show that
\begin{align*}
\max_{X_w(b)\neq\emptyset} d_w(b) - \dim X_w(b) = \begin{cases}k(k-1),&n=2k,\\
k^2,&n=2k+1.\end{cases}
\end{align*}

Since the dimension of ADLV is of general interest, we establish such a lower bound in the most general case, i.e.\ we no longer specialize to $\bG = \SL_n$.

\subsection{General setup}
Let $F$ be a non-archimedean local field with residue field $\BF_q$ and let $\breve F$ be the completion of the maximal unramified extension of $F$. We write $\Gamma$ for $\Gal(\overline F/F)$,  $\Gamma_0$ for the inertia subgroup of $\Gamma$ and $\sigma\in \Gamma$ for the Frobenius. Let $\bG$ be a connected reductive group over $F$ and $\breve G=\bG(\breve F)$. Let $\s$ be the Frobenius morphism on $\breve G$. We fix a $\s$-stable Iwahori subgroup $\breve I$ of $\breve G$. 
Let $Fl=\breve G/\breve I$ be the \emph{affine flag variety}. Let $\tW$ be the Iwahori-Weyl group of $\breve G$. Then we have a natural identification $\breve I \backslash \breve G/\breve I \cong \tW$ and the $\s$-action on $\breve G$ induces a natural action on $\tW$, which we still denote by $\s$. The extended affine Weyl group $\tW$ is the semidirect product of the finite Weyl group $W_0$ and the $\Gamma_0$-coinvariants of the cocharacter lattice $X_\ast(T)_{\Gamma_0}$.

For any $b \in \breve G$ and $w \in \tW$, we define the corresponding \emph{affine Deligne--Lusztig variety} in the affine flag variety $$X_w(b)=\{g \breve I \in \breve G/\breve I; g \i b \s(g) \in \breve I \dot w \breve I\} \subset Fl.$$ It is known that the affine Deligne--Lusztig variety $X_w(b)$ is a (probably empty) locally closed (perfect) scheme of locally finite type over the residue field of $\breve F$.  It is a general fact that one may reduce all questions regarding the geometry of affine Deligne--Lusztig varieties to the case where the group $\bG$ is quasi-split and of adjoint type \cite[Section~2]{Goertz2015}. Hence we assume from now on that $\bG$ satisfies these assumptions. In particular, this means that the finite Weyl group $W_0$ is stable under the action of $\sigma$.

We denote the Borovoi fundamental group of $\bG$ to be $\pi_1(\bG) = X_\ast(T)/\mathbb Z\Phi^\vee$, where $\Phi^\vee$ is the set of coroots. Then the Kottwitz point of $[b]\in B(\bG)$ is denoted $\kappa(b) \in \pi_1(\bG)_{\Gamma}$, which characterizes the connected components of the affine flag variety up to $\sigma$-action.

Write $w = xt^\mu y$ with $x,y\in W_0, \mu\in X_\ast(T)_{\Gamma_0}$ and $t^\mu y \in {}^{\BS} \tW$. Then we put $\eta_\sigma(w) = \sigma^{-1}(y)x$. Let $\nu_b\in X_\ast(T)_{\Gamma_0}\otimes\mathbb Q$ denote the dominant Newton point, $\de(b)\in \mathbb Z_{\ge 0}$ the defect and $2\rho\in X^\ast(T)^{\Gamma}$ the sum of the positive roots. Then we define the \emph{virtual dimension} of the pair $(w,[b])$ to be
\begin{align*}
d_w(b) = \frac 12\Bigl(\ell(w) + \ell(\eta_\sigma(w))-\langle \nu_b,2\rho\rangle - \de(b)\Bigr).
\end{align*}
Here we write $\ell$ for the length function of $\tW$ and $W_0$, and write $\ell_R$ for the \emph{reflection length} on $W_0$. The reflection length an element $u\in W_0$ as the smallest number $n$ such that $u$ is a product of $n$ reflections (not necessarily simple) in $W_0$. It is denoted
\begin{align*}
\ell_R(u) = \min\{n\mid \exists \alpha_1,\dotsc,\alpha_n\in \Phi \text{ such that } u = s_{\alpha_1}\cdots s_{\alpha_n}\}.
\end{align*}
Write $\mathcal O = \{w^{-1} w_0 \sigma(w)\mid w\in W_0\}$ for the $\sigma$-conjugacy class of the longest element $w_0$ and $\ell_R(\mathcal O) = \min\{\ell_R(u)\mid u\in\mathcal O\}$. For a complete description of $\ell_R(\mathcal O)$, we refer to \cite{He2021_dimension}.

Set $B(\mathbf{G})_w = \{ [b]\in B(\mathbf{G}); X_w(b)\ne \emptyset\}$ and let $[b]\in B(\bG)_w$. We know that $\dim X_w(b)\le d_w(b)$, and the goal of this section is to give a bound of the difference $d_w(b) - \dim X_w(b)$.

\begin{theorem}\label{thm:dimVirtDimBound}The maximum of the difference between virtual dimension and dimension, for all pairs $(w,[b])$ such that the corresponding affine Deligne-Lusztig variety is non-empty, is precisely given by
\begin{align*}
\max_{\substack{w\in \tW\\\relax[b]\in B(\bG)_w}} d_w(b) - \dim X_w(b) = \frac {\ell(w_0) - \ell_R(\mathcal O)}2.
\end{align*}
\end{theorem}

We summarize our proof strategy as follows. There is a unique maximal element in $B(\mathbf{G})_w$, denoted by $[b_{w}]$. It is called the \emph{generic} $\sigma$-conjugacy class of $w$, since it is the unique $\sigma$-conjugacy class such that the intersection $[b_w]\cap \breve Iw\breve I$ is dense in $\breve Iw\breve I$.
The existence of $[b_w]$ and a useful combinatorial characterization are obtained by Viehmann in \cite[Corollary~5.6]{Viehmann2014}.

By a deep result in arithmetic geometry, we will see that the difference of virtual dimension and dimension reaches its maximum, over $B(\bG)_w$, at $[b] = [b_w]$.

The advantage of working with $[b_w]$ is that we have an explicit formula for the dimension of $X_w(b_w)$. Combined with a description of $[b_w]$ via a certain combinatorial model, we then compute the difference $d_w(b_w)-\dim X_w(b_w)$. Finally, we re-write that upper bound in terms of the length and reflection length functions of the finite Weyl group $W_0$.

For $\bG = \SL_n$, one evaluates $\ell(w_0) = n(n-1)/2$  and $\ell_R(w_0) = \lfloor(n-1)/2\rfloor$ in order to obtain the upper bound as stated in the beginning of this section.
\subsection{Step 1: A purity result}
We denote the usual dominance order on $B(\bG)$ by $\le$. For $[b],[b']\in B(\bG)$, this means that $[b]\le [b']$ if and only if the Kottwitz points $\kappa(b),\kappa(b')$ agree as elements of $\pi_1(G)_{\Gamma}$ and the difference of Newton points $\nu_{b'}-\nu_b$ is a $\mathbb Q_{\ge 0}$-linear combination of positive coroots. Geometrically, this means that the subset $[b]\subset \breve G$ lies inside the closure of $[b']\subset \breve G$. 

One now may study increasing chains $[b] = [b_1]<[b_2]<\cdots<[b_{n+1}] = [b']$ in $B(\bG)$. By the work of Chai \cite[Theorem~7.4]{Chai2000}, we know that all maximal chains have the same length, given by
$$\text{length}([b],[b']) = \frac 12\Bigl(\langle \nu_{b'}-\langle\nu_b,2\rho\rangle +\de(b') - \de(b)\Bigr).$$

One may similarly ask for geometric properties of $\sigma$-conjugacy classes not inside $\breve G$, but inside the smaller subset $\breve Iw\breve I$. The intersection $[b]\cap \breve Iw\breve I$ is infinite-dimensional, but it is admissible in the sense of \cite{He2015}, so there is a well-defined notion of codimension in $\breve Iw\breve I$ and closure inside $\breve Iw\breve I$. There is a notion of relative dimension of $[b]\cap \breve Iw\breve I$, allowing us to express the dimension of $X_w(b)$ in terms of $\dim [b]\cap \breve Iw\breve I$.

For $[b]\in B(\bG)_w$, the closure of the Newton stratum $[b]\cap \breve Iw\breve I$ is contained in the union of all Newton strata $[b']\cap \breve Iw\breve I$ for $[b]\ge [b']\in B(\bG)_w$, but one cannot in general expect to have an equality.

To compare the dimensions of different Newton strata inside $\breve Iw\breve I$, we use the purity theorem.
This is a deep result in arithmetic geometry, developed by many experts, including de Jong, Viehmann and Hamacher. We won't recall the statement nor the proof, due to the level of technicalities involved, and instead refer to the discussion in \cite{Viehmann2020}.

The statement we use here is the due to Viehmann \cite[Lemma~5.12]{Viehmann2020}. It states that the codimension of $[b]\cap \breve Iw \breve I$ inside $\breve Iw\breve I$ is at most $\text{length}([b], [b_w])$. By \cite[Theorem~2.23]{He2015}, we know that this codimension is equal to $\dim X_w(b_w) - \dim X_w(b) + \<\nu_{b_w}-\nu_b,2\rho\>$. Thus
\begin{align}
\dim X_w(b_w) - \dim X_w(b) \le& \text{length}([b],[b_w]) -\<\nu_{b_w}-\nu_b,2\rho\>
\notag\\=&\frac 12\Bigl(\langle \nu_b - \nu_{b_w},2\rho\rangle +\de(b) - \de(b_w)\Bigr)
\notag\\=&d_w(b_w) - d_w(b).\label{eq:purityConsequence}
\end{align}
Hence $d_w(b) - \dim X_w(b) \le d_w(b_w) - \dim X_w(b_w)$.
In other words, the function $$B(\bG)_w\rightarrow\mathbb Z, [b]\mapsto d_w(b)-\dim X_w(b)$$ reaches its maximum at $[b] = [b_w]$. We now focus on this special case.

\subsection{Step 2: The quantum Bruhat graph}
Let $w\in \tW$.  Recall that $[b_w]\in B(G)_w$ denotes the generic $\sigma$-conjugacy class associated with $w$.
In this step, we calculate the difference $d_w(b_w) - \dim X_w(b_w)$ for arbitrary elements $w\in \tW$.

It is known that $\dim X_w(b_w) = \ell(w) - \langle \nu_{b_w},2\rho\rangle$, cf.\ \cite[Theorem~2.23]{He2015}. Thus we compute
\begin{align*}
d_w(b_w) - \dim X_w(b_w) = \frac 12\Bigl(-\ell(w)+\ell(\eta_\sigma(w)) +\langle \nu_{b_w},2\rho\rangle -\de(b_w)\Bigr).
\end{align*}
By \cite[Proposition~3.9]{Schremmer2022}, we know that $$\langle \nu_b,2\rho\rangle - \de(b)=\langle \lfloor b\rfloor,2\rho\rangle,$$ where $\lfloor b\rfloor\in X_\ast(T)_{\Gamma}$ is the best integral approximation of $[b]$ in the sense of \cite{Hamacher2018}. More specifically, $\lfloor b \rfloor$ is the unique element in $X_*(T)_{\G}$ such that 
\begin{itemize}
\item $\k(b) = \k(\lfloor b \rfloor)$ in $\pi_1(\bG)_{\G}$ and
\item $0\le \<\nu_b - \lfloor\nu_b\rfloor, \omega_o\> < 1$ for any $o\in \BS/\<\s\>$. Here, $\omega_o\in \mathbb Q\Phi$ is the unique weight whose pairing with a simple root $\alpha$ is given by $1$ if $s_\alpha\in o$ and $0$ otherwise.
\end{itemize}
It remains to compute this approximation $\lfloor b_w\rfloor$ for arbitrary elements $w\in \tW$. This is a result of Schremmer \cite[Theorem~4.2]{Schremmer2022}, generalizing earlier results which compute this quantity in special cases.

In order to understand these generic $\sigma$-conjugacy class $[b_w]$, the tool of choice for Schremmer's result and its predecessors is a finite combinatorial object associated with $\bG$, known as \emph{quantum Bruhat graph}. This graph was originally introduced by Brenti, Fomin and Postnikov \cite{Brenti1998} as a consequence of certain solutions to Yang-Baxter equations. While originally intended to study quantum cohomology, especially the quantum Chevalley-Monk formula, it has since been found useful in a number of contexts, such as Kirillov-Reshetikhin crystals \cite{Lenart2015} and Bruhat order of affine Weyl groups \cite{Lam2010}. The calculation of the generic $\sigma$-conjugacy class of affine Weyl group elements is related to the Bruhat order via a result of Viehmann \cite[Corollary~5.6]{Viehmann2014}. The resulting connection between the quantum Bruhat graph and the generic $\sigma$-conjugacy class of sufficiently regular affine Weyl group elements was first discovered by Mili\'cevi\'c \cite{Milicevic2021}.

The quantum Bruhat graph is defined as follows: By definition, $\text{QBG}(\Phi)$ is a directed graph, whose set of vertices is given by the finite set $W_0$. Its edges are of the form $w\rightarrow ws_\a$ for $w\in W_0$ and $\a\in \Phi^+$ whenever one of the following conditions is satisfied:
\begin{itemize}
\item $\ell(ws_\a) = \ell(w)+1$ or
\item $\ell(ws_\a) = \ell(w)+1-\langle \a^\vee,2\rho\rangle$.
\end{itemize}
Edges satisfying the first condition are called \emph{Bruhat} edges, whereas edges satisfying the second condition are called \emph{quantum} edges (hence the graph's name). It is common to draw the graph with the vertical position of the vertices corresponding to the length, with the longest element on top and the identity element at the bottom. Then the Bruhat edges go upwards whereas the quantum edges go downwards. This is the quantum Bruhat graph of type $A_2$:

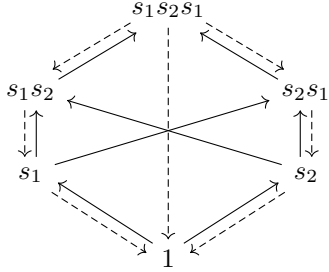
\def\qecolor{dashed}
\def\seshift{0.5ex}
\def\quantumEdge{dashed}
\def\shortEdgeShiftRight{0.5ex}
\begin{align*}
\begin{tikzcd}[ampersand replacement=\&,column sep=2em]
\&s_1 s_2 s_1\ar[ddd,\qecolor]
\ar[dl,\qecolor,shift right=\seshift]\ar[dr,\qecolor,shift left=\seshift]\\
s_1 s_2\ar[ur,shift right=\seshift]\ar[d,\qecolor,shift right=\seshift]\&\&
s_2 s_1\ar[ul,shift left=\seshift]\ar[d,\qecolor,shift left=\seshift]\\
s_1\ar[u,shift right=\seshift]\ar[urr]\ar[dr,\qecolor,shift right=\seshift]\&\&s_2\ar[u,shift left=\seshift]\ar[ull]\ar[dl,\qecolor,shift left=\seshift]\\
\&1\ar[ru,shift left=\seshift]\ar[lu,shift right=\seshift]
\end{tikzcd}
\end{align*}
It is known that the quantum Bruhat graph is (strongly) connected. Hence we may write $d_{\text{QBG}}(u,v)\in\mathbb Z_{\ge 0}$ for the length of a shortest path from $u$ to $v$ in the quantum Bruhat graph, where $u,v\in W_0$.

The description of the generic $\sigma$-conjugacy class $[b_w]$ in terms of $w$ uses the quantum Bruhat graph and the decomposition $w = xt^\mu y$ as above. The latter decomposition is not canonical if $w$ is not very regular, so we might have to vary it slightly. This can be done using the notion of length positive elements, as introduced by Schremmer \cite{Schremmer2022}.

Let $w = t^\lambda z\in\tW$. We say that $v\in W_0$ is \emph{length positive} with respect to $w$ if all positive roots $\alpha\in\Phi^+$ satisfy
$$\<z^{-1}\lambda,v\alpha\> + \delta(zv\alpha)-\delta(v\alpha) \ge 0. $$
Denote the set of length positive elements by $\LP(w)\subseteq W_0$. If we write $w = xt^\mu y$ with $t^\mu y\in {}^{\BS} \tW$ as above, then $y^{-1}$ is always length positive, i.e.\ $y^{-1}\in \LP(w)$.

With this setup, we can summarize the main result of \cite{Schremmer2022} as follows:
\begin{theorem}{{\cite[Theorem 4.2 and Lemma 4.4]{Schremmer2022}}}\label{Felix}
Let $w =t^\lambda z\in \tW$. Let $[b_w]$ the maximal element in $B(\bG)_w$. Then  

\begin{align*}& \< \lfloor b_w\rfloor , 2\rho\> = \ell(w) -\min_{v\in \LP(w)}d_{\text{QBG}}\bigl(v,\s(zv)\bigr).\pushQED{\qed}\qedhere\popQED\end{align*}
\end{theorem}
In view of the above calculation, we obtain
\begin{align}
d_w(b_w) - \dim X_w(b_w) = \frac 12\Bigl(\ell(y \s(x)) - \min_{v\in \LP(w)}d_{\text{QBG}}\bigl(v,\s(zv)\bigr)\Bigr).\label{eq:genericUpperBound}
\end{align}
A priori, since there are infinitely many elements in $\tW$, it is not clear whether or not the left-hand side of the above equation has an upper bound. However, since $W_0$ is a finite group, the right hand-side (and thus also the left-hand side) of equation~\eqref{eq:genericUpperBound} has an upper bound.

\subsection{Step 3: The reflection length as an upper bound}
From equations \eqref{eq:genericUpperBound} and \eqref{eq:purityConsequence}, we see that
\begin{align*}
\max_{[b]\in B(\bG)_w} d_w(b) - \dim X_w(b) = \frac 12\Bigl(\ell(y \sigma(x)) - \min_{v\in \LP(w)}d_{\text{QBG}}\bigl(v,\s(zv)\bigr)\Bigr).
\end{align*}
It remains to compute the maximum of this expression over all $w\in\tW$. A major difficulty in explicitly evaluating the right-hand side of the above expression is the condition $v\in \LP(w)$. In this section, we relax this condition to $v\in W_0$, thus obtaining the upper bound
\begin{align*}
\max_{[b]\in B(\bG)_w} d_w(b) - \dim X_w(b)\le \frac 12\max_{x,y,v\in W} \Bigl(\ell(y \sigma(x)) - d_{\text{QBG}}\bigl(v,\s(xyv)\bigr)\Bigr).
\end{align*}

From the definition of the quantum Bruhat graph, we see
\begin{align*}
\min_{v\in W_0}d_{\text{QBG}}\bigl(v,\s(xyv)\bigr) \ge
\min_{v\in W_0} \ell_R(v^{-1} \sigma(xyv)).
\end{align*}
We summarize
\begin{align*}
\max_{\substack{w\in \tW\\\relax[b]\in B(\bG)_w}}
d_w(b)-\dim X_w(b) \le &\max_{v,x,y\in W_0} \frac 12\Bigl(\ell(y\sigma(x)) - \ell_R(v^{-1}\sigma(xyv))\Bigr)
\intertext{Writing $u=yv$, we can re-write this as}\cdots=&\max_{u,x,y\in W_0} \frac 12\Bigl(\ell(y\sigma(x)) - \ell_R(u^{-1} y\sigma(x)\sigma(u))\Bigr)
\\=&\max_{u,\eta\in W_0}\frac 12\Bigl(\ell(\eta) - \ell_R(u^{-1} \eta\sigma(u))\Bigr).
\intertext{
If $\eta\neq w_0$, we find a simple reflection $s$ with $\ell(\eta s) = \ell(\eta)+1$. Then certainly $\ell_R(u^{-1}\eta s\sigma(u))\le \ell_R(u^{-1} \eta\sigma(u))+1$. So when searching for the above maximum, we may replace $\eta$ by $\eta s$ until $\eta = w_0$. Thus, we can simplify the above expression to
}
\cdots =& \max_{u\in W_0}\frac 12\Bigl(\ell(w_0) - \ell_R(u^{-1} w_0 \sigma(u))\Bigr).
\end{align*}
We proved that
\begin{align*}
\max_{\substack{w\in \tW\\\relax[b]\in B(\bG)_w}}
d_w(b)-\dim X_w(b) = \frac 12\Bigl(\ell(y \sigma(x)) - \min_{v\in \LP(w)}d_{\text{QBG}}\bigl(v,\s(zv)\bigr)\Bigr)  \le \frac{\ell(w_0)-\ell_R(\mathcal O)}2.
\end{align*}
obtaining the upper bound claimed in Theorem~\ref{thm:dimVirtDimBound}.

The reader will find a peculiar similarity to the paper \cite{He2021_dimension} of He and Yu. They study a similar maximization problem in \cite[Lemma~4.3, Theorem~5.1]{He2021_dimension}, proving that
\begin{align*}\max_{x,y\in W_0} \Bigl(\ell(\sigma^{-1}(y)x) - d_{\mathrm{QBG}}(x,y^{-1})\Bigr) = \ell(w_0) - \ell_R(\mathcal O).
\end{align*}

\subsection{Step 4: Explicit construction of the lower bound}
We saw in the previous step that 
\begin{align*}
\max_{\substack{w\in \tW\\\relax[b]\in B(\bG)_w}}
d_w(b)-\dim X_w(b) = \frac 12\Bigl(\ell(y \sigma(x)) - \min_{v\in \LP(w)}d_{\text{QBG}}\bigl(v,\s(zv)\bigr)\Bigr)  \le \frac{\ell(w_0)-\ell_R(\mathcal O)}2.
\end{align*}
In this section, we prove that equality holds, by explicitly constructing an element $w\in \tW$ and $v\in \LP(w)$ such that $\eta_\sigma(w) = y\sigma(x) = w_0$ and $d_{\text{QBG}}\bigl(v,\s(xyv)\bigr) = \ell_R(\mathcal O)$. We do this construction in a case-by-case fashion, first considering the case where $\bG$ is quasi-simple over $\breve F$, meaning that the root system is connected.

\subsubsection{The split and $\sigma = \Ad(w_0)$ cases}\label{subs:splitAdw0}
Consider first the case where the action of $\sigma$ on $\Phi$ is either the identity map (i.e.\ $\bG$ is residually split) or equal to the action of $-w_0$. In either case, we obtain $\ell_R(\mathcal O) = \ell_R(w_0)$. From \cite[section~5]{Sadhukhan2021}, we know that $\ell_R(w_0) = d_{\text{QBG}}(w_0, 1)$. Then $w := t^{2\rho^\vee} w_0$ satisfies $\eta_\sigma(w) = w_0$ and $\LP(w) = \{w_0\}$, from which one obtains $d_w(b_w) - \dim X_w(b_w) = \frac{\ell(w_0) - \ell_R(\mathcal O)}2$, completing the proof of the theorem.

Following the classification of root systems, one sees that $\sigma$ is given by one of the above choices \emph{unless} $\bG$ is non-split of type $D_n$ with $n$ even.
\subsubsection{The case ${}^2 D_{2k}$.}\label{subs:2D2k}Suppose that $\bG$ is of type $D_{2k}$ with $k\ge 2$ and the image of $\sigma$ in $\Aut(\Phi)$ has order $2$. Label the simple roots as $\alpha_1,\dotsc,\alpha_{2k}$ such $\alpha_i$ is connected to $\alpha_{i+1}$ in the Dynkin diagram of $D_{2k}$ for all $i=1,\dotsc,2k-2$. Then $\sigma$ interchanges the roots $\alpha_{2k}$ and $\alpha_{2k-1}$, while fixing all other roots. The element $w_0\in W_0$ is central. One computes $\ell_R(\mathcal O) = 2k-2$ \cite[section~5.7]{He2021_dimension}.

Define $x = s_{2k-1}, y = w_0 s_{2k}, v = w_0\in W_0$ and $\mu \in X_\ast(T)_{\Gamma_0}$ as $\mu = 2\rho_K^\vee$, the sum of all positive coroots of the sub-root system spanned by $K = \{\alpha_1,\dotsc,\alpha_{2k-1}\}$. Then $w = xt^{\mu} y \in \tW$ satisfies $y = y^{-1} \in \LP(w)$ and $\ell(w, y \alpha_{2k}) = 0$, hence $v\in \LP(w)$. We get $\eta_\sigma(w) = y\sigma(x) = w_0$.

It suffices to show that $d_{\text{QBG}}(v, \sigma(xyv)) \le \ell_R(\mathcal O) = 2k-2$. For this, we compute

\begin{align*}
d_{\text{QBG}}(v, \sigma(xyv)) &= d_{\text{QBG}}(w_0, s_{2k-1} s_{2k}) = d_{\text{QBG}}(s_{2k-1} w_0, s_{2k})\\&=d_{\text{QBG}}(s_{2k} s_{2k-1} w_0, 1),
\end{align*}
where we applied \cite[Lemma~7.7]{Lenart2015} twice. Since $w_0, s_{2k-1}$ and $s_{2k}$ centralize each other, we write $s_{2k} s_{2k-1} w_0 = w_0 s_{2k-1} s_{2k}$.

Denote the longest root of $\Phi^+$ by $\theta$. Applying \cite[Lemma~7.7]{Lenart2015} again, we conclude
\begin{align*}
d_{\text{QBG}}(w_0 s_{2k-1} s_{2k},1) = 1+d_{\text{QBG}}(s_\theta w_0 s_{2k-1} s_{2k},1).
\end{align*}
Let $J = \{\alpha_1,\alpha_3,\dotsc,\alpha_{2k}\}$ such that the longest element of $W_J$ is equal to $s_\theta w_0$. Then $s_\theta w_0 s_{2k-1} s_{2k} = s_1 w_0(J') s_{2k-1} s_{2k}$ where $J' = J\setminus\{\alpha_1\}$. We conclude
\begin{align*}
d_{\text{QBG}}(w_0 s_{2k-1} s_{2k},1) = 2+d_{\text{QBG}}(w_0(J') s_{2k-1} s_{2k},1).
\end{align*}
If $k=2$, one checks that $w_0(J')s_{2k-1} s_{2k}=1$ so that the desired identity follows. Otherwise, we have $k\geq 3$ and observe that $J'$ defines a $\sigma$-stable subroot system of type $D_{2k-2}$. In an inductive sense, we may assume that $d_{\text{QBG}}(w_0(J')s_{2k-1} s_{2k},1)=2(k-1)-2$ has already been established. Then also $d_{\text{QBG}}(w_0 s_{2k-1} s_{2k},1)=2k-2$ follows immediately.

This completes the proof that for the $w$ constructed above, we have $d_w(b_w) - \dim X_w(b) = \frac {\ell(w_0) -  \ell_R(\mathcal O)}2$, establishing Theorem~\ref{thm:dimVirtDimBound} for groups of type $D_{2k}$ where $\sigma$ has order $2$.

\subsubsection{The case ${}^3 D_4$}\label{subs:3D4}
If $\bG$ is of type $D_4$ and $\sigma$ has order $3$, enumerate the simple roots $\alpha_1,\dotsc,\alpha_4$ such that $\sigma(\alpha_1) = \alpha_3$, $\sigma(\alpha_3) = \alpha_4$ and $\sigma(\alpha_4) = \alpha_1$. Then one chooses $x = s_4, y = w_0 \sigma^{-1}(x) = w_0 s_1$ and $v = w_0$ as above. Moreover, we choose $\mu=2\rho_J^\vee$ with $J = \{\alpha_2,\alpha_3,\alpha_4\}$. Then one chooses $w = x t^\mu y\in \tW$ and $v\in \LP(w)$ to get the same conclusion as in the example above.

\subsubsection{The general case}
Following \cite[section~5.3]{He2021_dimension}, we can reduce to the case where $\bG$ is quasi-simple over $F$. In other words, we assume that $\s$ acts transitively on the set of irreducible components of $W_0$. We have $W_0 = W_0' \times \ldots W_0'$ with $\ell$ irreducible components. There is a length-preserving group automorphism $\s'$ on $W_0'$ with $\s(w_1,\ldots, w_{\ell}) = (\s'(w_{\ell}),w_1,\ldots,w_{\ell-1} )$. Let $w_0'$ be the longest element of $W_0'$ and let $\CO'$ be the $\s'$-conjugacy class of $w_0'$ in $W_0'$. We distinguish the $\ell$ even case and the $\ell$ odd case as in \cite[Section~5.4]{He2021_dimension}.

If $\ell$ is even. Let $u = (1,w_0',1,w_0',\ldots,1,w_0')$, then $w_0 = u\s(u)$ and hence $\ell_R(\CO) = 0$. Let $x = y=(w_0',1,w_0',1,\ldots,w_0',1)$ and $\mu = (2\rho^{\vee},0,2\rho^{\vee},0,\ldots,2\rho^{\vee},0)$. Consider $w = x t^{\mu} y$. Then $y\s(x) = w_0$, $v = y^{-1}u \in \LP(w)$ and $d_{\text{QBG}}(y^{-1}u, \s(xu)) = d_{\text{QBG}}(w_0,w_0) = 0$. Therefore, $w$ satisfies the condition.

If $\ell$ is odd. Then by \cite[Section~5.4.2]{He2021_dimension}, $\ell_R(\CO) = \ell_R(\CO')$. By the result proved above in subsections \ref{subs:splitAdw0}, \ref{subs:2D2k} and \ref{subs:3D4}, one can find $w' = x't^{\mu'}w_0'\s(x'^{-1})$ and $u'\in W_0'$ such that $\s(x')w_0'u'\in \LP(w')$ and $d_{\text{QBG}}(\s(x')w_0'u',\s'(x'u')) = \ell_R(\CO')$. Now let 
\begin{align*}
u &= (u',1,1,\ldots,1),\\
v &= ( \s'(x')w_0' u',x'w_0',x' ,x'w_0',x' ,\ldots,x'w_0',x'),\\
x &= (x',x'w_0',x',x'w_0',\ldots,x',x'w_0',x'),\\
y &=(w_0'\s'(x'^{-1}),w_0'x'^{-1},x'^{-1},w_0'x'^{-1},x'^{-1},\ldots,w_0'x'^{-1},x'^{-1}),\\   
\mu &= (\mu',\rho_{J_1}^{\vee},\rho_{J_2}^{\vee},\rho_{J_1}^{\vee},\rho_{J_2}^{\vee},\ldots,\rho_{J_1}^{\vee},\rho_{J_2}^{\vee} ),
\end{align*}
where $J_1 = \{i;x'w_0'(\a_i)<0\}$, $J_2 = \{i;x'(\a_i)<0\}$ and $\rho_{J_1}^{\vee}$ and $\rho_{J_2}^{\vee}$ are sum of fundamental coweights corresponding to $J_1$ and $J_2$ respectively. Consider $w = x t^{\mu} y$. Then $y\s(x) = w_0$, $y^{-1}u \in \LP(w)$ and $$d_{\text{QBG}}(y^{-1}u, \s(xu)) = d_{\text{QBG}}(\s(x')w_0'u',\s'(x'u')) =  \ell_R(\CO')=\ell_R(\CO).$$ Therefore, $w$ satisfies the condition. Theorem~\ref{thm:dimVirtDimBound} is fully proved.

\subsection{Final comments}
This concludes the theoretical discussion of the difference $d_w(b) - \dim X_w(b)$.
It is noteworthy that, even though the dimension of affine
Deligne--Lusztig varieties can be fully determined by the combinatorial
algorithm presented in section~\ref{sec:ADLVAlgo}, we had to employ
tools from quantum cohomology and algebraic geometry to prove
Theorem~\ref{thm:dimVirtDimBound}. This situation is typical in research
on affine Deligne--Lusztig varieties and illustrates why the field
encompasses more than just analyzing a single algorithm on affine Weyl
groups.

\section{Conclusion}
\label{sec:conclusion}

We used machine learning to study a central unsolved problem in pure mathematics, namely the geometry of affine Deligne--Lusztig varieties. In this section, we want to discuss the potential of this new research method and share some practical insights on the interdisciplinary collaboration.

Our project required an interdisciplinary research group, consisting of experts in machine learning and specialized researchers of the mathematical problem at hand. This joint expertise allowed us to find machine learning models by using subject--specific knowledge, as well as interpret their behavior from the perspective of a subject matter expert. After exchanging explanations on the machine learning models used and the mathematical problem to be studied, we established a common understanding of the material. We could even delve into highly technical questions about modeling specific discrete functions, such as the dimension of affine Deligne-Lusztig varieties, using neural networks.

While the research method employed in this project can not lead to new mathematical \emph{proofs}, it still offers new insights, raises intriguing questions, and leads to conjectures. Starting with a naive view of the mathematical problem, we were able to rediscover some of the most crucial tools and invariants developed by the mathematical community with substantial amount of time and effort. This demonstrates that this pipeline has the potential to accelerate research. 

Additionally, we identified new avenues of research that could be explored using established mathematical methods. By analyzing the linear model developed in \S 4, we formulated a new conjecture that has since been proven in \S 6.

We would like to emphasize some important requirements that contributed to our success, as a recommendation to other mathematical researchers (even from very different fields) who might benefit from this new research method:

\begin{itemize}
\item \textbf{Selection of Problems}: Choose problems that are easy to generate data for but difficult to find patterns in. In our case, the mathematical problem we studied was computable for a large number of examples. While the general goal of ``better understanding the geometry of affine Deligne-Lusztig varieties" is too vague, we were able to define specific numerical invariants, such as non-emptiness, dimension, and the number of irreducible components. With a substantial amount of computed data points ranging from thousands to millions, machine learning becomes a powerful tool to study the patterns.

\item \textbf{Machine Learning Model Selection}: The choice of the machine learning model is crucial. While large neural networks may offer higher accuracy, they can reduce interpretability and make identifying patterns challenging. Complex function forms can also introduce noise in sensitivity analysis. Hence, it is important to strike a balance between model complexity and interpretability. Furthermore, the choice of loss functions and training algorithms may also play a vital role.

\item \textbf{Leveraging Prior Knowledge}: Prior knowledge plays a significant role in the success of this approach. It aids in selecting appropriate features and models and helps researchers explain and evaluate the results. Understanding the problem itself is crucial for distinguishing meaningful patterns from noise, especially when facing counterintuitive results.

\item \textbf{Moderate Technical Requirements}: The technical requirements are relatively low, and researchers do not need an overly complex machine learning setup. In our experiments, we found that even with millions of data points, these small networks could acquire well-trained models in just a few minutes when utilizing a single GPU. Training times on conventional laptops were generally less than an hour. However, the generation of data can be time-consuming. Despite not having an effective parallel implementation, generating millions of data points on a single CPU may still take several days. 

\item \textbf{Flexibility in Concerned Function}: One of the advantages of machine learning models is their ability to handle functions that are not smooth, continuous, or even discrete. In our case, the problem did not satisfy some common assumptions in machine learning, such as input feature independence. While this introduced occasional challenges, we could often find a way to alleviate it, squeezing out useful information about the machine learning model and the underlying mathematical problem in the process. 

\item \textbf{Effective Communication between Groups}: Interdisciplinary communication, particularly between groups steeped in pure mathematics and machine learning domains, is of paramount importance. In our research, we found that in-depth discussions on data generation, regularization and fidelity terms significantly enhanced our comprehension of the numerical outcomes.

While individual experiments can be executed swiftly, the subsequent analysis and interpretation of results often demand substantial time. The process inherently necessitates periodic revisions and discussions, leading to a cyclical pattern of repeated experimentation. Remarkably, the analysis and interpretation phase tends to consume the lion's share of time in this iterative process.

Enhancing the mutual understanding between these expert groups could significantly expedite this process. A deeper appreciation of each other's domains can foster clearer communication, leading to more efficient analysis and potentially quicker identification of necessary experimental modifications. This cross-disciplinary understanding, therefore, serves as a catalyst for accelerating the overall research process.
Conversely, the interdisciplinary research process itself improves the mutual understanding naturally over time.
\end{itemize}


\printbibliography

\end{document}